\documentstyle[12pt,amstex]{article}
\title{Besov Spaces of Self-affine Lattice Tilings and Pointwise regularity}
\author{ Koichi Saka }
\date{\empty}
\newtheorem{Theorem}{\bf Theorem}
\newtheorem{Lemma}{\bf Lemma}
\newtheorem{Proposition}{\bf Proposition}
\newtheorem{Corollary}{\bf Corollary}
\begin{document}
\maketitle
\bigskip
\begin{abstract}
We investigate Besov spaces of self-affine tilings of ${\Bbb R}^{n}$ and 
discuss various characterizations of those Besov spaces. 
We see what is a finite set of functions which generates the Besov spaces 
from a view of multiresolution approximation on self-affine lattice tilings 
of ${\Bbb R}^{n}$. 
Using this result we give a generalization of already known 
characterizations of Besov spaces 
given by wavelet expansion
 and we apply to study the pointwise H${\ddot {\rm o}}$lder 
space. Furthermore  we give descriptions of scaling exponents measured 
by Besov spaces, and estimations of a pointwise H${\ddot {\rm o}}$lder 
exponent to compute the pointwise scaling exponent of 
several oscillatory functions.



\end{abstract}

\section{Introduction}
There are many ways to characterize Besov spaces. 
Among them in  the discrete version are regular wavelet expansion, 
Littlewood-Paley decomposition, polynomial approximation, spline 
approximation, mean oscillation, and difference operator 
(See  [9], [13] and [15]).
 In this paper we give these characterizations in context of 
 self-affine lattice tilings of ${\Bbb R}^{n}$ 
  and we apply to study these pointwise versions. 
 In particular we see to give most of these characterizations 
 in a framework of multiresolution approximation 
 on self-affine lattice tilings of ${\Bbb R}^{n}$. 
 We also give conditions of finitely many functions which 
 generate the Besov spaces of self-affine lattice tilings of 
 ${\Bbb R}^{n}$ in a view of multiresolution approximation scheme. 
 This result is a generalization of characterizations of Besov spaces 
 given by regular wavelet functions and by spline functions.(See  [3], [12] and [15]). 
 Moreover we apply to give descriptions of scaling exponents by 
 characterizations of the Besov space, and we also consider 
 a pointwise  H${\ddot {\rm o}}$lder exponent 
 of oscillatory functions given by 
 a multiresolution approximation series in self-affine lattice tilings 
 of ${\Bbb R}^{n}$.

  The plan of sections in our paper is as follows:

  In the second section 
we introduce self-affine lattice tilings of ${\Bbb R}^{n}$ which  
arise in many contexts, particularly, in fractal geometry 
 and in construction of wavelet bases. See [14] for a survey on related 
 topics. 
We define  Besov spaces of self-affine lattice tilings, 
and give its characterizations and its pointwise versions.

In the third section  
we consider a multiresolution analysis $\{ V_{l} \}$ generated 
 by finitely many  functions associated 
 with a self-affine lattice tiling. 
 We give properties of Besov space norms defined by approximation errors 
 associated with $\{ V_{l} \}$.

 In the fourth section 
we  give some conditions  
of finitely many functions which characterize the Besov space by 
multiresolution approximation on self-affine lattice tilings 
of ${\Bbb R}^{n}$. 
 We apply this result to give a generalization of 
 characterizations of Besov spaces given by regular wavelet functions 
 and  by spline functions,  
and we also give characterizations of 
the  pointwise H${\ddot {\rm o}}$lder space 
by multiresolution approximation.

In the fifth section 
we give descriptions of scaling exponents of 
global and poitwise regularity 
by characterizations of the Besov space.   
We give properties of a pointwise H${\ddot {\rm o}}$lder exponent for 
a multiresolution approximation series in self-affine lattice tilings 
and apply to compute a pointwise H${\ddot {\rm o}}$lder exponent 
of several oscillatory functions.

 We use  $C$ to denote a positive constant different in each occasion. 
 But it will depend on the parameter appearing in each problem. 
 The same notations $C$ are not necessarily 
 the same on any two occurrences. 

\section{Self-affine lattice tilings and Besov spaces}

\noindent

Let $\Gamma$ be a lattice in ${\Bbb R}^{n}$,
 that is, $\Gamma$ is an image of the integer lattice ${\Bbb Z}^{n}$ 
under some nonsingular linear transformation 
and let $M$ be a dilation matrix, that is, all eigenvalues of $M$ have 
absolute values greater than one and $M$ preserves the lattice $\Gamma$: 
$\ \ M\Gamma \subset \Gamma$. This implies that $|\det M|=m$ is a positive integer 
greater than one and $m$ is the order of the quotient space $\Gamma/M\Gamma.$ 
We say that a compact set $T$ generates a self-affine tiling $\{ T+\gamma \}
_{\gamma \in \Gamma}$ if  
\begin{eqnarray}
\cup_{\gamma \in \Gamma} (T+\gamma) = {\Bbb R}^{n} \ \ {\rm disjoint} \ \ 
{\rm a.e.}\, \nonumber \\ 
\cup_{\gamma \in \Gamma_{0}}(T+\gamma)=MT \ \ {\rm disjoint}\ \ {\rm a.e.} 
\end{eqnarray}
where  $\Gamma_{0}$ is a finite subset of $\Gamma$ consisting of 
representatives for disjoint cosets in $\Gamma/M\Gamma$. The set $\Gamma_{0}$ 
is called a set of digits and the compact set $T$ is called a self-affine tile.
 The self-affine tile $T$ has nonempty interior $T^{o}$. 
 For convenience we restrict our attention to 
 the case $\Gamma = {\Bbb Z}^{n}$. 
 In this case the dilation matrix $M$ has integer entries.

 For $1 \leq p \leq \infty$, let ${\cal L}^{p} = {\cal L}^{p}({\Bbb R}^{n})$ 
 be the linear space of all functions $\phi$ for which 

\begin{eqnarray}
|\phi|_{p}=(\int_{T}(\sum _{\nu \in {\Bbb Z}^{n}}|\phi(x-\nu)|)^{p}dx)^{1/p}
 < \infty.
 \end{eqnarray}
with the usual modification for $p= \infty$. 
 Clearly, 
 ${\cal L}^{p} \subset L^{p}({\Bbb R}^{n})$ and 
 ${\cal L}^{\infty} \subset {\cal L}^{p} \subset {\cal L}^{q} \subset 
 {\cal L}^{1} = L^{1}({\Bbb R}^{n})$ for $1 \leq q \leq p \leq \infty$.
If $\phi \in L^{p}({\Bbb R}^{n})\ (1 \leq p \leq \infty)$ is compactly supported, then $\phi \in {\cal L}^{p}$. Furthermore, we observe that if there are 
constants $C>0$ and $\delta >0$ such that 
$
|\phi (x)| \leq C(1+|x|)^{-n-\delta} 
$ for all $x \in {\Bbb R}^{n}$ 
then $\phi \in {\cal L}^{\infty}$.    

A finite subset $\Phi =\{ \phi_{1}, \ldots , 
\phi_{N} \}$ of ${\cal L}^{\infty}$ is said to have $L^{p}$-stable shifts 
$(1 \leq p \leq \infty)$, if there are constants $C_{1}>0$ and 
$C_{2}>0$ such that for any 
sequences $c_{j} \in l^{p}({\Bbb Z}^{n})\ \ (j=1, \ldots , N)$, 
\begin{eqnarray*}
C_{1}\sum^{N}_{j=1}||c_{j}||_{l^{p}} 
\leq ||\sum^{N}_{j=1}\sum_{\nu \in {\Bbb Z}^{n}}c_{j}(\nu)\phi_{j}(x-\nu)||_{p} \leq C_{2}\sum^{N}_{j=1}||c_{j}||_{l^{p}}.
\end{eqnarray*}
From now  those equivalences shall be described as 
\begin{eqnarray*}
\sum^{N}_{j=1}||c_{j}||_{l^{p}} 
\sim ||\sum^{N}_{j=1}\sum_{\nu \in {\Bbb Z}^{n}}c_{j}(\nu)\phi_{j}(x-\nu)||_{p} .
\end{eqnarray*}

\bigskip

\noindent
${\bf Theorem  \ A}$ ([6]). 
{\it
For a finite subset $\Phi=\{ \phi_{1}, \ldots, \phi_{N} \}$ of 
${\cal L}^{\infty}$, we have following equivalent conditions:

\noindent
{\rm (i)}
 $\Phi$ has $L^{2}$-stable shifts,

\noindent
{\rm (ii)}
 $\Phi$ has $L^{p}$-stable shifts for $1 \leq p \leq \infty$,

\noindent
{\rm (iii)}
 there is a set of functions $\tilde{\Phi}=\{ \tilde{\phi}_{1}, 
 \ldots, \tilde{\phi}_{N} \}$
 in ${\cal L}^{\infty}$, dual to $\Phi$ in the sense that 
 \begin{eqnarray*}
 \int \phi_{j}(x-\mu)\bar{\tilde{\phi}}_{k}(x-\nu) dx = 
 \delta_{\mu\nu}\delta_{jk}, \ \ j,k=1, \ldots, N, \ \ 
 \mu,\nu \in {\Bbb Z}^{n},
 \end{eqnarray*}
where $\delta$ is the Kronecker's symbol.
}

\bigskip

Let $\Pi=\{ T+\nu \}_{\nu \in {\Bbb Z}^{n}}$ be a self-affine lattice 
tiling of ${\Bbb R}^{n}$ with a dilation matrix $M$.  
For a nonnegative integer $k$, we denote the 
function $p_{\alpha}$ with $|\alpha| \leq  k,\alpha \in {\Bbb Z}^{n}_{+}$, 
where ${\Bbb Z}_{+}$ is the set of all nonnegative integers, given by 
\begin{eqnarray}
p_{\alpha}(x) =x^{\alpha}, \ \  x \in T^{o} \nonumber \\ 
p_{\alpha}(x)=0  \ \ {\rm otherwise}.
\end{eqnarray}
Since $\Phi=\{ p_{\alpha}\}_{|\alpha|\leq k}$ of ${\cal L}^{\infty}$ 
has $L^{2}$-stable shifts, there is a set of functions $\tilde{\Phi}=
\{ \tilde{p}_{\alpha} \}_{|\alpha|\leq k}$ dual to $\Phi$.

Let $Q_{0}$ be a translate of the tile $T$ containing the origin 
as an interior point and let $p^{\prime}_{\alpha}$, 
$\tilde{p}^{\prime}_{\alpha}$ be corresponding translates of 
$p_{\alpha}$, $\tilde{p}_{\alpha}$ respectively.
For $Q_{l}(x_{0})=M^{-l}Q_{0}+x_{0}$, we write 
\[
p^{Q_{l}(x_{0})}_{\alpha}(x)=m^{l/2}p^{\prime}_{\alpha}(M^{l}(x-x_{0})), \ \ \ 
\tilde{p}^{Q_{l}(x_{0})}_{\alpha}(x)=m^{l/2}\tilde{p}^{\prime}_{\alpha}
(M^{l}(x-x_{0})) 
\]
\begin{eqnarray}
P_{Q_{l}(x_{0})}f(x)=\sum_{|\alpha|\leq k}
\langle f,\tilde{p}^{Q_{l}(x_{0})}_{\alpha}\rangle p^{Q_{l}(x_{0})}_{\alpha}(x).
\end{eqnarray}

\bigskip

We define 
\begin{eqnarray}
{\rm osc}^{k}_{p}f(x,l)=\inf_{P\in {\Bbb P}^{k}}(\frac{1}{|Q_{l}(x)|}
\int_{Q_{l}(x)}|f(y)-P(y)|^{p}dy)^{1/p}
\end{eqnarray}
and 
\begin{eqnarray*}
\Omega^{k}_{p}f(x,l)=(\frac{1}{|Q_{l}(x)|}\int_{Q_{l}(x)}
|f(y)-P_{Q_{l}(x)}f(y)|^{p}dy)^{1/p}
\end{eqnarray*}
where $Q_{l}(x)=M^{-l}Q_{0}+x$ and $P_{Q_{l}(x)}f$ is given in (4), 
and $|Q_{l}(x)|$ is the volume element of $Q_{l}(x)$, and 
${\Bbb P}^{k}$ is the linear space of all polynomials of degree no greater 
than $k$ on ${\Bbb R}^{n}$.

\bigskip

\noindent
{\bf Definition}. 
Let $\lambda_{0}$ be the least value of absolute values of eigenvalues of the 
dilation matrix $M$.
 Given  $s>0$,  $k$  a nonnegative integer with $k+1 >s$ 
 and $1 \leq p,q \leq \infty$. A function $f$ is said to 
belong to the Besov space $B^{s}_{pq}(M)$ if 
\begin{eqnarray}
||f||_{B^{s}_{pq}(M)}=
||f||_{p}+(\sum_{l=0}^{\infty}(\lambda_{0}^{ls}
||{\rm osc}_{p}^{k}f(\cdot,l)||_{p})^{q})^{1/q} 
<\infty.
\end{eqnarray}
with the usual modification for $q=\infty$. 
We note that the above definition is independent of the choice of  
 nonnegative integers $k$ with $k+1 >s$ and ${\rm osc}_{p}^{k}$ in the 
 definition can be replaced by ${\rm osc}_{1}^{k}$. 
 We can see $W^{p}_{k+1}({\Bbb R}^{n}) \subset B^{s}_{pq}(M)$ if $s<k+1$.
When the dilation matrix $M$ is $\lambda_{0}$-times of the identy $Id$ with 
$\lambda_{0} > 1$, 
the above Besov space coincides the usual Besov space on ${\Bbb R}^{n}$.

\bigskip

\noindent
{\rm\bf Remark 1}. 
We have the embedding theorem : $B^{\beta}_{p \xi}(M) \subset 
B^{\alpha}_{p \eta}(M)$ for $\beta > \alpha >0$, $1 \leq \xi, \eta \leq 
\infty$ and $1 \leq p \leq \infty$, and $B^{\alpha}_{p \xi}(M) \subset 
B^{\alpha}_{p \eta}(M)$ for $\alpha > 0$, $1 \leq \xi \leq \eta \leq \infty$ 
and $1 \leq p \leq \infty$.

\bigskip

Let $\triangle_{u}f$ denote the difference operator 
$\triangle_{u}f(x)= f(x+u)-f(x)$.
Let us choose positive constants $r$ and $d$  such that 
\begin{eqnarray}
\{ u \in {\Bbb R}^{n}: |u| < r \} \subset Q_{0} 
\subset \{ u \in {\Bbb R}^{n}: |u| < dr \}.
\end{eqnarray}

\bigskip

\begin{Theorem}. 
Given  $s>0$ , a nonnegative integer $k$ with $k+1> s$ and 
$1 \leq p,q \leq \infty$,
we have equivalent ones of the Besov space norm given in {\rm (6)},
if one of them exists, with the usual modification for $q=\infty$, 
\begin{eqnarray*}
  \lefteqn{
||f||_{B^{s}_{pq}(M)}
} \\
 &\sim& 
||f||_{p}+(\sum_{l=0}^{\infty}(\lambda_{0}^{ls}
||\Omega_{p}^{k}f(\cdot,l)||_{p})^{q})^{1/q} 
\equiv |||f|||_{1}, \\
 &\sim&
||f||_{p}+(\sum^{\infty}_{l=0}(\lambda_{0}^{ls}
\sup_{(k+1)|M^{l}u|<r/2}||\triangle^{k+1}_{u}f||_{p}
)^{q})^{1/q}
\equiv |||f|||_{2}.
\end{eqnarray*}
\end{Theorem}

{\rm\bf Proof.} 
Since ${\rm osc}_{p}^{k}f(x, l) \sim \Omega_{p}^{k}f(x, l)$ 
from [2], 
the equivalence of $||f||_{B^{s}_{pq}(M)}$ and $|||f|||_{1}$ is obvious.

We shall prove for any $f$ such that $|||f|||_{2} <\infty $, 
\begin{eqnarray*}
||\Omega_{p}^{k}f(\cdot, l)||_{p} \leq C(
\sup_{(k+1)|M^{l}u|<r/2}||\triangle^{k+1}_{u}f||_{p}
+\lambda_{0}^{-l(k+1)}||f||_{p}).
\end{eqnarray*}

 We choose a function $\chi$ in $C^{\infty}_{c}({\Bbb R}^{n})$ such that 
 $\int |\chi(u)| du =1$ and supp$\chi \subset \{ u \in {\Bbb R}^{n}: \ \ 
 |u| < r/2(k+1) \}$. We write 
 $\chi_{l}(u)=m^{l}\chi(M^{l}u),\ \ h_{l}(x)=\int(f(x)-\triangle_{u}^{k+1}f(x))
 \chi_{l}(u)du$ . 
  Then we have  
 \begin{eqnarray*}
\lefteqn{\Omega_{p}^{k}f(x, l) \leq (\displaystyle{\frac{1}{|Q_{l}(x)|}}
\int_{Q_{l}(x)}|f(y)-h_{l}(y)|^{p} dy)^{1/p} + 
(\displaystyle{\frac{1}{|Q_{l}(x)|}}
\int_{Q_{l}(x)}|h_{l}(y)-P_{Q_{l}(x)}h_{l}(y)|^{p} dy)^{1/p}}  \\ 
&+& (\displaystyle{\frac{1}{|Q_{l}(x)|}}\int_{Q_{l}(x)}
|P_{Q_{l}(x)}h_{l}(y)-P_{Q_{l}(x)}f(y)|^{p} dy)^{1/p}  \\ 
&\leq&  C((\displaystyle{\frac{1}{|Q_{l}(x)|}}
\int_{Q_{l}(x)}|f(y)-h_{l}(y)|^{p} dy)^{1/p} + \Omega_{p}^{k}h_{l}(x ,l)) 
 \equiv C(I_{1}(x)+I_{2}(x)).
 \end{eqnarray*}
We have :
 \[
 ||I_{1}||_{p} \leq ||(m^{l}\int_{M^{-l}Q_{0}}
 (\int |\triangle^{k+1}_{u}f(\cdot+y)|
 |\chi_{l}(u)| du)^{p} dy)^{1/p}||_{p} \leq 
 C\sup_{(k+1)|M^{l}u|< r/2}||\triangle^{k+1}_{u}f||_{p}. 
 \]

Let $q_{z}$ be the $k$-th Taylor polynomial of $h_{l}$ about 
$z \in {\Bbb R}^{n}$.
To estimate  $I_{2}$, we use 
\[
|\partial^{\beta}h_{l}(x)|  
\leq 
C\sum_{e=1}^{k+1}
\int_{|u|< r/2(k+1)}|f(x-eM^{-l}u)|du 
 \leq
  C\int_{|u|< r/2}|f(x-M^{-l}u)|du.
\]
Hence we get an estimate:
\begin{eqnarray*}
\lefteqn{
{\rm osc}_{p}^{k}h_{l}(x, l) \leq C(m^{l}\int_{Q_{l}(x)}
|h_{l}(y)-q_{x}(y)|^{p} dy)^{1/p} }\\ 
&\leq&  C(m^{l}\int_{Q_{l}(x)}|\int_{0}^{1}\sum_{|\beta|=k+1}
\displaystyle{\frac{k+1}{\beta !}}\partial^{\beta}
h_{l}(x+t(y-x))(1-t)^{k}(y-x)^{\beta}dt|^{p}dy)^{1/p} \\
&\leq& 
C(m^{l}\int_{M^{-l}Q_{0}}(\int_{0}^{1}\int_{|u|<r/2}|f(x+ty-M^{-l}u)||y|^{k+1}
dt du)^{p} dy)^{1/p}. 
 \end{eqnarray*}
Hence, we get an estimate : 
 \begin{eqnarray*}
||I_{2}||_{p}  
\leq 
C(m^{l}\int_{M^{-l}Q_{0}}|y|^{(k+1)p}||f||_{p}^{p} dy)^{1/p} 
 \leq 
 C||f||_{p}\lambda_{0}^{-l(k+1)}.
\end{eqnarray*}
Now we combine the estimates of $I_{1}$ and $I_{2}$ to write 
\begin{eqnarray*}
||\Omega_{p}^{k}f(\cdot, l)||_{p} \leq C(||I_{1}||_{p}+||I_{2}||_{p}) 
\leq C(\sup_{(k+1)|M^{l}u|< r/2}||\triangle^{k+1}_{u}f||_{p}+
\lambda_{0}^{-l(k+1)}||f||_{p}).
\end{eqnarray*}
This shows   $|||f|||_{1} \leq C|||f|||_{2}$.

We shall show the converse.
We can show from [5] for $|M^{l}u|<r/2(k+1)$,
\[
|\triangle_{u}^{k+1}f(x)| 
= |\triangle_{u}^{k+1}(f-P_{Q_{l}(x)}f)(x)|
 \leq C\sum_{e=0}^{k+1}\sum_{j=l}^{\infty}\Omega_{p}^{k}f(x+eu,j).
\]
Hence we have 
\[
\sup_{(k+1)|M^{l}u| < r/2}||\triangle_{u}^{k+1}f||_{p} 
\leq 
C\sum_{j=l}^{\infty}||\Omega_{p}^{k}f(\cdot, j)||_{p} 
\]
and we get the estimate 
\[
|||f|||_{2} \leq 
C|||f|||_{1}.
\]
This completes the proof of Theorem 1.

\bigskip

If $0 < s< k+1$ for a nonnegative integer $k$ and $1 \leq p,q \leq \infty$, 
then for $x \in {\Bbb R}^{n}$, a function $f \in T^{s}_{pq}(x)$ means that 
\begin{eqnarray*}
(\sum_{l=0}^{\infty}(\lambda_{0}^{ls}{\rm osc}^{k}_{p}f(x,l))^{q})^{1/q}
< \infty
\end{eqnarray*}
with the usual modification for $q=\infty$.

\bigskip

\noindent
{\rm\bf Remark 2}. 
We have the embedding theorem : $T^{\beta}_{p \xi}(x) \subset 
T^{\alpha}_{p \eta}(x)$ for $\beta > \alpha >0$, $1 \leq \xi, \eta \leq 
\infty$ and $1 \leq p \leq \infty$, and 
$T^{\alpha}_{p \eta}(x) \subset T^{\alpha}_{p \xi}(x), \ 
T^{\alpha}_{\xi q}(x) \subset T^{\alpha}_{\eta q}(x)$ 
for $\alpha > 0$, $1 \leq \eta \leq \xi \leq \infty$ 
and $1 \leq p, q \leq \infty$.

\bigskip

We have a poinwise version of Theorem 1, 
which is proved by the same way as the proof of Theorem 1.

\bigskip

\noindent
{\bf Corollary. \ }{\it  
Given  $s>0$ , a nonnegative integer $k$ with $k+1> s$ and 
$1 \leq p,q \leq \infty$. Then 
for $x \in {\Bbb R}^{n}$ following properties of 
a bounded function f are equivalent,  
 with the usual modification for $q=\infty$, 
\begin{eqnarray*}
&{\rm (i)}& 
f \in T^{s}_{pq}(x), \\
 &{\rm (ii)}& 
(\sum_{l=0}^{\infty}(\lambda_{0}^{ls}
\Omega^{k}_{p}f(x,l))^{q})^{1/q} < \infty, \\
 &{\rm (iii)}&
(\sum^{\infty}_{l=0}(\lambda_{0}^{ls}
\sup_{(k+1)|M^{l}u|< r/2}(\displaystyle{\frac{1}{|Q_{l}(x)|}}
\int_{Q_{l}(x)}|\triangle^{k+1}_{u}f(y)|^{p}dy)^{1/p})^{q})^{1/q}<\infty.
\end{eqnarray*}
}

We will define the Littlewood-Paley decomposition. 
Let us $\lambda_{0} > 1$ and $\varphi$ a function in the Schwartz class 
${\cal S}({\Bbb R}^{n})$ 
with the following properties: 
supp $\hat{\varphi} \subset \{ \xi \in {\Bbb R}^{n} :\ |\xi| \leq 1 \}$ 
and $\hat{\varphi}(\xi) =1$  on 
$\{ \xi \in {\Bbb R}^{n} : \ |\xi|\leq \lambda_{0}^{-1} \}$. 
Let $\psi(x) = \lambda_{0}^{n}\varphi(\lambda_{0}x)-\varphi(x)$. 
Let $\varphi_{l}(x) = \lambda_{0}^{ln}\varphi(\lambda_{0}^{l}x)$, 
$S_{l}f=f*\varphi_{l}$, 
$\psi_{l}(x)=\lambda_{0}^{ln}\psi(\lambda_{0}^{l}x)$ and 
$f_{l}= f*\psi_{l}$ for $l=0,\ 1,\ 2,\ \ldots$. Then 
for $f \in {\cal S}'$ we have Littlewood-Paley decomposition:
\begin{eqnarray}
f=\varphi*f + \sum_{l=0}^{\infty}\psi_{l}*f \equiv S_{0}f +
\sum_{l=0}^{\infty}f_{l}.
\end{eqnarray}

\bigskip

\noindent
{\bf Theorem \ B} ([13]). {\it 
Suppose that a dilation matrix is of the form M = $\lambda_{0}$Id with 
$\lambda_{0} > 1$. Let $1 \leq p, q \leq \infty$ and $s > 0$. 
Then we have equivalence of norms if one of them exit, 
for Littlewood-Paley decomposition given in {\rm (8)}, 
with the usual modification $q= \infty$: 
\begin{eqnarray*}
&{\rm (i)}& \lefteqn{ \ ||f||_{B^{s}_{pq}(M)} ,} \\
\sim 
&{\rm (ii)}& 
 ||f||_{p} + 
(\sum_{l=0}^{\infty}(\lambda_{0}^{ls}||f-S_{l}f||_{p})^{q})^{1/q}, \\ 
\sim
&{\rm (iii)}&
 ||S_{0}f||_{p} +(\sum_{l=0}^{\infty}(\lambda^{ls}_{0}
||f_{l}||_{p})^{q})^{1/q}. 
\end{eqnarray*}
}

We write $T_{\infty \infty}^{s}(x)=C^{s}(x)$. 
The following statement is a pointwise version of Theorem B.

\bigskip

\begin{Proposition}. 
Suppose that a dilation matrix is of the form M = $\lambda_{0}$Id with 
$\lambda_{0} > 1$. Let $ s > 0$. 
Then for $x \in {\Bbb R}^{n}$, following properties of 
a bounded function f for 
Littlewood-Paley decomposition given in {\rm (8)} are equivalent: 
\begin{eqnarray*}
&{\rm (i)}&  f \in C^{s}(x) , \\
&{\rm (ii)}& 
|f(y)-S_{l}f(y)| \leq C(\lambda_{0}^{-l}+|x-y|)^{s}\ \  {\rm for}
 \ {\rm all}\ l \geq 0. 
\end{eqnarray*}
\end{Proposition}

{\bf Proof}. 
We will prove that (i) implies (ii). Let us $k+1 > s$. By [1: Theorem 2] 
 we have 
\[
f(x) - S_{l}f(x) = \int \triangle^{k+1}_{u}f(y)\varphi_{l}(u)du 
- \sum_{e=1}^{k}\left( \begin{array}{c}
k \\ e 
\end{array} \right)
(-1)^{e}f*\psi_{l}^{e}(y)
\]
where $\varphi^{e}(y)=e^{-n}\varphi (e^{-1}y)$ and 
$\psi^{e}=\varphi^{e}-\varphi^{e+1}$. 
From [1: Lemma 2] there exist functions 
$\psi^{e1}, \ldots , \psi^{en}$ in ${\cal S}({\Bbb R}^{n})$ such that 
\[
\psi_{l}^{e}(x)=\sum_{i=1}^{n}\triangle^{k+1}_{c\lambda^{-l}_{0}e_{i}}
\psi^{ei}_{l}(x)
\]
where $c=\lambda_{0}^{-1}e\pi$ and $e_{1}, \ldots , e_{n}$ 
are the canonical basis vectors in ${\Bbb R}^{n}$ 
to write 
\[
f*\psi_{l}^{e}=\sum_{i=1}^{n}f*(\triangle^{k+1}_{c\lambda_{0}^{-l}e_{i}}
\psi_{l}^{ei})
=\sum_{i=1}^{n}(\triangle^{k+1}_{c\lambda_{0}^{-l}e_{i}}f)*\psi^{ei}_{l}.
\]
We get fom the corollary of Theorem 1
\begin{eqnarray*}
\lefteqn{|f*\psi^{e}_{l}(y)| \leq C\sum_{i=1}^{n}
\int |\triangle^{k+1}_{c\lambda_{0}^{-l}e_{i}}f(y-z)\psi^{ei}_{l}(z)|dz 
} \\ 
&\leq&  
C\sum_{i=1}^{n}\int(|x-y+z|+c\lambda_{0}^{-l})^{s}
|\psi^{ei}_{l}(z)|dz 
\leq 
C\sum_{i=1}^{n} \int(|x-y|+|z|+\lambda_{0}^{-l})^{s}|\psi^{ei}_{l}(z)|dz 
\\ 
&\leq& 
C\sum_{i=1}^{n}(\lambda_{0}^{-l}+|x-y|)^{s}\int(1+|z|)^{s}|\psi^{ei}(z)|dz 
\leq  C(\lambda_{0}^{-l}+|x-y|)^{s}. 
\end{eqnarray*}
Furthermore we see 
\begin{eqnarray*}
\lefteqn{\int |\triangle^{k+1}_{u}f(y)\varphi_{l}(u)|du 
\leq C\int (|x-y|+|u|)^{s}|\varphi_{l}(u)|du 
} \\ 
&\leq& 
C(\lambda_{0}^{-l}+|x-y|)^{s}\int(1+|u|)^{s}|\varphi(u)|du 
\leq C(\lambda_{0}^{-l}+|x-y|)^{s}
\end{eqnarray*}
 These give the estimate
 \begin{eqnarray*}
\lefteqn{ |f(y)-S_{l}f(y)| 
} \\
&\leq& C\int|\triangle_{u}^{k+1}f(y)\varphi_{l}(u)|du + C\sum_{e=1}^{k}|f*\psi_{l}^{e}(y)|
 \leq C(\lambda_{0}^{-l}+|x-y|)^{s}.
\end{eqnarray*} 
 This completes the implication of (ii) from (i).

We will show that (ii) implies (i). 
We choose a positive integer $l_{1}$ such that $\lambda_{0}^{l_{1}} > d+1$ 
where $d$ is given in (7). 

We have 
\[
|\triangle^{k+1}_{u}f(y)| \leq |\triangle^{k+1}_{u}(f-S_{l}f)(y)| + 
|\triangle^{k+1}_{u}S_{l}f(y)|
\]
We will give an estimate of $|\triangle^{k+1}_{u}(f-S_{l}f)(y)|$ 
to see 
\begin{eqnarray*}
\lefteqn{\sup_{Q_{l+l_{1}}(x)}\sup_{(k+1)|M^{l+l_{1}}u|< r/2}
|\triangle^{k+1}_{u}(f-S_{l}f)(y)| 
} \\ 
&\leq& C\sup_{Q_{l+l_{1}}(x)}\sup_{(k+1)|M^{l+l_{1}}u|<r/2}
\sum_{e=0}^{k+1}|(f-S_{l}f)(y+eu)| 
\leq C\sup_{Q_{l}(x)}|(f-S_{l}f)(y)|.
\end{eqnarray*}
We will give an estimate of $|\triangle^{k+1}_{u}S_{l}f(y)|$ to write 
\[
|\triangle^{k+1}_{u}S_{l}f(y)| \leq C\sum_{|\beta|=k+1}|u|^{k+1}
\int_{0}^{1}\cdots \int_{0}^{1}
|\partial^{\beta}S_{l}f(y+(\theta_{1}+\cdots +\theta_{k+1})u)|
d\theta_{1}\ldots d\theta_{k+1}.
\]
We have by (ii)
\[
|f_{l}(y)| \leq |f(y)-S_{l+1}f(y)|+|f(y)-S_{l}f(y)|
\leq C(\lambda_{0}^{-l}+|x-y|)^{s}.
\]
Hence, to estimate $|\partial^{\beta}S_{l}f(y)|$ with $|\beta|=k+1$ 
we use, by Bernstein's inequality, 
\begin{eqnarray*}
\lefteqn{|\partial^{\beta}S_{l}f(y)|\leq \sum_{j=0}^{l-1}
|\partial^{\beta}f_{j}(y)| +| \partial^{\beta}S_{0}f(y)| 
} \\ 
&\leq& C(\sum_{j=0}^{l-1}\lambda_{0}^{(k+1)j}(\lambda_{0}^{-j}+|x-y|)^{s}
+||f||_{\infty}) \leq C\lambda_{0}^{(k+1)l}(\lambda_{0}^{-l}+|x-y|)^{s}.
\end{eqnarray*}
From those we have 
\[
\sup_{(k+1)|M^{l}u| <r/2}|\triangle_{u}^{k+1}S_{l}f(y)| \leq 
C\sup_{(k+1)|M^{l}u|<r/2}|u|^{k+1}\lambda_{0}^{(k+1)l}(\lambda_{0}^{-l}+|x-y| 
+|u|)^{s} \leq C(\lambda_{0}^{-l}+|x-y|)^{s}.
\]
Hence we get an estimate 
\begin{eqnarray*}
\lefteqn{\sup_{Q_{l+l_{1}}(x)}\sup_{(k+1)|M^{l+l_{1}}u|<r/2}|
\triangle_{u}^{k+1}f(y)|
} \\
&\leq& C\sup_{Q_{l}(x)}|(f-S_{l}f)(y)|  
+ \sup_{Q_{l}(x)}\sup_{(k+1)|M^{l}u|<r/2}|\triangle_{u}^{k+1}S_{l}f(y)|
\leq C\lambda_{0}^{-ls}.
\end{eqnarray*}
These complete the proof of Proposition 1 by the corollary of Theorem 1.

\bigskip

A following corollary may be proved by the same way as 
in the proof of Proposition 1.

\bigskip

\noindent
{\bf Corollary}. {\it Suppose that a dilation matrix $M=\lambda_{0}Id$. 
Let f be a bounded function. 
If $f \in C^{s}(x)$, then it holds

\bigskip

${\rm (iii)}\ \ \ \ \ |f_{l}(y)| \leq C(\lambda_{0}^{-l}+|x-y|)^{s}\ \ $ 
 for all $l \ \geq \  0$.

\bigskip

Conversely, if it holds for $s > s^{\prime} > 0$,

\bigskip

${\rm (iii)^{\prime}}\ \ \ \ \ |f_{l}(y)| \leq C\lambda_{0}^{-ls}
(1+\lambda_{0}^{l}|x-y|)^{s^{\prime}}\ \ $  for all $l \ \geq \  0$,

\bigskip

\noindent
then $f \in C^{s}(x)$.
}

\section{Multiresolution approximation}

\noindent
Let $\Pi$ denote a self-affine lattice tiling 
$\{ T+\nu \}_{\nu \in {\Bbb Z}^{n}}$ with a dilation matrix $M$. 
For an integer $l$ and a finite subset $\Phi = \{ \phi_{1}, 
\ldots, \phi_{N} \}$ of ${\cal L}^{\infty}$ with $L^{2}$-stable shifts, 
we define operators $P_{l}f$ given by 
\begin{eqnarray}
P_{l}f(x) = \sum^{N}_{j=1}\sum_{\nu \in {\Bbb Z}^{n}}
m^{l}\langle f, {\tilde{\phi}}_{j}(M^{l}\cdot-\nu) \rangle 
\phi_{j}(M^{l}x-\nu) 
\end{eqnarray}
where $\langle f,{\tilde{\phi}}_{j}(M^{l}\cdot-\nu)\rangle =
\int f(y)\bar{\tilde{\phi}}_{j}(M^{l}y-\nu) \ dy$
 and $\tilde{\Phi} =\{ \tilde{\phi}_{1}, \ldots, \tilde{\phi}_{N} \}$ is 
 dual to $\Phi$ in Theorem A.

Let $V_{0}^{p}=\{ \sum_{j=1}^{N}\sum_{\nu \in {\Bbb Z}^{n}}a_{j}(\nu)
\phi_{j}(x-\nu): \ \ a_{j} \in l^{p}({\Bbb Z}^{n}) \}$  
and let $V_{l}^{p}=\{ f(M^{l}x) : f \in V_{0}^{p} \}$. 
Then for $1 \leq p  \leq \infty$, the operator $P_{l}$ is a bounded 
projection operator of $L^{p}({\Bbb R}^{n})$ onto 
 $V_{l}^{p} \ \  (1 \leq p \leq \infty)$ in the sense that $P_{l}f =f $
 for any $f \in V_{l}^{p}$. 
 We say $\Phi =\{ \phi_{1}, \ldots, \phi_{N} \}$ of ${\cal L}^{\infty}$ 
 is $M$-refinable if there exist sequences 
 $c_{jk} \in l^{1}({\Bbb Z}^{n})  \ \ (1 \leq j,k \leq N)$ such that 
 \[
 \phi_{j}(x)= \sum^{N}_{k=1}\sum_{\nu \in {\Bbb Z}^{n}} c_{jk}(\nu)
 \phi_{k}(Mx-\nu), \ \ \  x \in {\Bbb R}^{n} ,  \ \ \ j=1,\ldots, N.
 \]

A following theorem implies that $\{ V_{l}^{p} \}$ is a multiresolution 
analysis in $L^{p}({\Bbb R}^{n})$ for $1 \leq p < \infty$.

\bigskip

\noindent
{\bf Theorem C} ( [6] and [16]). 
{\it
  If a finite subset $ \Phi$ of ${\cal L}^{\infty}$ is $M$-refinable and 
has $L^{2}$-stable shifts, then the sequence of sets 
$\{ V_{l}^{p} \}\ (1 \leq p \leq \infty)$ satisfies following properties:

\noindent
{\rm (i)}
 $f \in V_{0}^{p} \Leftrightarrow f(x-\nu) \in V_{0}^{p}$ 
for all $\nu \in {\Bbb Z}^{n}$ ,

\noindent
{\rm (ii)}
 $f \in V_{l}^{p} \Leftrightarrow f(Mx) \in V_{l+1}^{p}$,

\noindent
{\rm (iii)}
 $ \cdots \subset V_{l}^{p} \subset V_{l+1}^{p} \subset \cdots$ ,

\noindent
{\rm (iv)}
$\cap_{l \in {\Bbb Z}}V_{l}^{p}=\{ 0 \}
\ \ (1 \leq p < \infty)$,

\noindent
{\rm (v)}
$ \cup^{\infty}_{l=0}V_{l}^{p}$ is dense in $L^{p}({\Bbb R}^{n}) 
\ \ (1 \leq p < \infty).$
}

\bigskip

Given a function $f$ in $L^{p}({\Bbb R}^{n})\ \  (1 \leq p \leq \infty)$, 
$\sigma^{p}_{l}(f)$ denotes the error of $L^{p}$-approximation from 
$V_{l}^{p}$ in $L^{p}({\Bbb R}^{n})$:
\begin{eqnarray}
\sigma^{p}_{l}(f) = \inf \{ ||f-S||_{p} : S \in V_{l}^{p} \}.
\end{eqnarray}
 
\noindent
Clearly we have the following equivalence:
\begin{eqnarray*}
\sigma^{p}_{l}(f) \sim ||f-P_{l}f||_{p},\ \ f \in L^{p}({\Bbb R}^{n})\ \ 
(1 \leq p \leq \infty ).
\end{eqnarray*}
 Given  $s>0$,  $\lambda > 1$  and $1 \leq p,q \leq \infty$. 
 A function $f$ is said to 
belong to  $B^{s, \lambda}_{pq}(\Phi)$ if 
\begin{eqnarray}
||f||_{B^{s, \lambda}_{pq}(\Phi)}=||f||_{p}
+(\sum^{\infty}_{l=0}(\lambda^{ls}\sigma^{p}_{l}(f))^{q})^{1/q}<\infty
\end{eqnarray}
with the usual modification when  $q=\infty$.

Let 
\begin{eqnarray}
R_{l}f = P_{l+1}f - P_{l}f, \ \ l=0, 1, \ldots . 
\end{eqnarray}
We put 
\begin{eqnarray}
P_{0}f(x)=\sum_{j=1}^{N}\sum_{\nu \in {\Bbb Z}^{n}}
a_{j0}(\nu)\phi_{j}(x-\nu), \ \ \ 
R_{l}f(x)=\sum_{j=1}^{N}\sum_{\nu \in {\Bbb Z}^{n}}
a_{j(l+1)}(\nu)\phi_{j}(M^{l+1}x-\nu). 
\end{eqnarray}
Since $\Phi$ has stable shifts, we have 
\begin{eqnarray}
||P_{0}f||_{p} \sim \sum_{j=1}^{N}||a_{j0}||_{l^{p}},\ \ \ 
||R_{l}f||_{p} \sim m^{-(l+1)/p}\sum_{j=1}^{N}||a_{j(l+1)}||_{l^{p}},
\ \  l=0,1,\ldots.
\end{eqnarray}  
Then for $f \in B_{pq}^{s, \lambda}(\Phi)$ we have 
\[
f(x) =P_{0}f(x) + \sum_{l=0}^{\infty}R_{l}f(x) \equiv 
\sum_{j=1}^{N}\sum_{l=0}^{\infty}
\sum_{\nu \in {\Bbb Z}^{n}}a_{jl}(\nu)\phi_{j}(M^{l}x-\nu).
\]
Moreover from [15,  Theorem 5.10] there exists an associated set of wavelets 
$\{ \psi^{\epsilon}_{j} \}_{j=1,\ldots,N}^{\epsilon =1,\ldots,m-1}$, 
that is, 
$\{ \psi^{\epsilon}_{j}(x-\nu) \}_{j=1,\ldots,N, \nu \in {\Bbb Z}^{n}}
^{\epsilon =1,\ldots,m-1}$ is 
an orthonormal basis in $W_{0} = V_{1}^{2} \ominus V_{0}^{2}$ 
in $L^{2}({\Bbb R}^{n})$ , whose wavelet expansion of a function 
$f \in L^{2}({\Bbb R}^{n})$ is given by 
\begin{eqnarray}
f(x) = \sum_{j=1}^{N}\sum_{\nu \in {\Bbb Z}^{n}}a_{j0}(\nu)\phi_{j}(x-\nu) 
+\sum_{j=1}^{N}\sum_{\epsilon=1}^{m-1}\sum_{l=0}^{\infty}
\sum_{\nu \in {\Bbb Z}^{n}}b_{jl}^{\epsilon}(\nu)m^{l/2}\psi^{\epsilon}
_{j}(M^{l}x-\nu)
\end{eqnarray}
where 
\begin{eqnarray}
a_{j0}(\nu)=\langle f(y), \tilde{\phi}_{j}(y-\nu) \rangle ,\ \ 
b^{\epsilon}_{jl}(\nu)=\langle f(y), m^{l/2}\psi^{\epsilon}_{j}(M^{l}y-\nu) 
\rangle.
\end{eqnarray}
Then we have 
\[
P_{0}f(x)=\sum_{j=1}^{N}\sum_{\nu \in {\Bbb Z}^{n}}a_{j0}(\nu)\phi_{j}(x-\nu),
\] 
\begin{eqnarray*}
R_{l}f(x)=\sum_{j=1}^{N}\sum_{\epsilon=1}^{m-1}\sum_{\nu \in {\Bbb Z}^{n}}
b_{jl}^{\epsilon}(\nu)m^{l/2}\psi_{j}^{\epsilon}(M^{l}x-\nu), l=0,1,\ldots.
\end{eqnarray*}
When $m >(n+1)/2$, there exist  $\psi^{\epsilon}_{j} \in {\cal L}^{\infty}$ 
and 
\begin{eqnarray*}
||R_{l}f||_{p} \sim m^{l(1/2-1/p)}\sum_{j=1}^{N}\sum_{\epsilon=1}^{m-1}
||b_{jl}^{\epsilon}||_{l^{p}}\ \ (1 \leq p \leq \infty).
\end{eqnarray*}

\begin{Theorem}.  Assume that a finite subset 
$\Phi=\{ \phi_{1}, \ldots, \phi_{N} \}$ of ${\cal L}^{\infty}$ 
is $M$-refinable and has $L^{2}$-stable shifts. 
Given $\lambda > 1$ and $\alpha > 0$, there are equivalences of 
the norm $||f||_{B^{\alpha,\lambda}_{pq}(\Phi)}$ given in {\rm (11)}, 
if one of them exits, for any $1 \leq p ,q \leq \infty$,
with the usual modification for $q=\infty$:

\noindent
{\rm (i)} 
$ ||f||_{p} +(\sum_{l=0}^{\infty}(\lambda^{l\alpha}||f-P_{l}f||_{p})^{q})
^{1/q}
$,

\noindent
{\rm (ii)} 
$ ||P_{0}f||_{p} + (\sum_{l=0}^{\infty}(\lambda^{l\alpha}||R_{l}f||_{p})
^{q})^{1/q}
$,

\noindent
{\rm (iii)} $(\sum_{l=0}^{\infty}(\lambda^{l\alpha}
m^{-l/p}\sum_{j=1}^{N}||a_{jl}||_{l^{p}})^{q})^{1/q}$,

\noindent
where $\{ a_{jl} \}$  are given in {\rm (13)}.

\noindent
{\rm (iv)} $
\inf (\sum_{l=0}^{\infty}(\lambda^{l\alpha}m^{-l/p}\sum_{j=1}^{N}||c_{jl}||
_{l^{p}})^{q})^{1/q}$
 where the infimum is taken over all admissible $L^{p}$-convergent 
 representations 
\[
f(x) = \sum_{j=1}^{N}\sum_{l=0}^{\infty}\sum_{\nu \in {\Bbb Z}^{n}}c_{jl}(\nu)
\phi_{j}(M^{l}x-\nu),
\]

\noindent
{\rm (v)} $\sum_{j=1}^{N}||a_{j0}||_{l^{p}} +
(\sum_{l=0}^{\infty}(\lambda^{l\alpha}
m^{l(1/2-1/p)}\sum_{j=1}^{N}\sum_{\epsilon=1}^{m-1}||b_{jl}^{\epsilon}
||_{l^{p}})^{q})^{1/q}\ \ \ {\rm when}\ \  m > (n+1)/2$,

\noindent
where $\{ a_{j0} \}$ and $\{ b_{jl}^{\epsilon} \}$ are given in {\rm (16)}.
\end{Theorem}

{\bf Proof}.  
Those equivalences of (i), (ii), (iii), (v) and (11) can be proved 
from easy routine using Hardy's inequality. 
We will prove only the equivalence of (iv) and (11). 
We prove that 
\[
||f||_{p}+(\sum_{l=0}^{\infty}(\lambda^{l\alpha}\sigma^{p}_{l}(f))^{q})^{1/q}
\leq C(\sum_{l=0}^{\infty}(\lambda^{l\alpha}m^{-l/p}
\sum_{j=1}^{N}||c_{jl}||_{l^{p}})^{q})^{1/q} < \infty
\]
when 
$ f(x) = \sum_{j=1}^{N}\sum_{l=0}^{\infty}\sum_{\nu \in {\Bbb Z}^{n}}
c_{jl}(\nu)\phi_{j}(M^{l}x-\nu)$ 
is $L^{p}$-convergent.

We put $f(x)=\sum_{l=0}^{\infty}F_{l}(x), \ \ 
F_{l}(x) = \sum_{j=1}^{N}\sum_{\nu \in {\Bbb Z}^{n}}
c_{jl}(\nu)\phi_{j}(M^{l}x-\nu)$. Then we see 
$F_{l} \in V_{l}^{p}$ and 
$||F_{l}||_{p} \sim m^{-l/p}\sum_{j=1}^{N}||c_{jl}||_{l^{p}}$.

Hence we have 
\begin{eqnarray*}
\lefteqn{
\sigma^{p}_{l_{0}}(f) = \sigma^{p}_{l_{0}}(\sum_{l=0}^{\infty}F_{l})
} 
\\ 
&\leq& 
\sum_{l=0}^{\infty}\sigma_{l_{0}}^{p}(F_{l}) 
=\sum_{l=0}^{l_{0}}\sigma_{l_{0}}^{p}(F_{l}) + \sum_{l=l_{0}+1}^{\infty}
\sigma_{l_{0}}^{p}(F_{l}) = \sum_{l=l_{0}+1}^{\infty}\sigma_{l_{0}}^{p}(F_{l})
\leq \sum_{l=l_{0}+1}^{\infty}||F_{l}||_{p}.
\end{eqnarray*}
From Hardy's inequality this implies 
\begin{eqnarray*}
\lefteqn{
||f||_{p}+
(\sum_{l_{0}=0}^{\infty}(\lambda^{l_{0}\alpha}\sigma_{l_{0}}^{p}(f))^{q})^{1/q} \leq 
\sum_{l=0}^{\infty}||F_{l}||_{p}+
(\sum_{l_{0}=0}^{\infty}(\lambda^{l_{0}\alpha}\sum_{l=l_{0}+1}^{\infty}
||F_{l}||_{p})^{q})^{1/q} } \\ 
&\leq& 
C(\sum_{l=0}^{\infty}(\lambda^{l\alpha}||F_{l}||_{p})^{q})^{1/q} 
\leq C(\sum_{l=0}^{\infty}(\lambda^{l\alpha}m^{-l/p}
\sum_{j=1}^{N}||c_{jl}||_{l^{p}})^{q})^{1/q}. 
\end{eqnarray*}
From (iii), we have 
\begin{eqnarray*}
||f||_{B^{\alpha, \lambda}_{pq}(\Phi)} 
\sim
(\sum_{l=0}^{\infty}(\lambda^{l\alpha}m^{-l/p}\sum_{j=1}^{N}
||a_{jl}||_{p})^{q})^{1/q}.
\end{eqnarray*}
This implies the equivalence of (iv) and (11).

\bigskip

\noindent
\begin{Proposition}.   
Given $k+1 >s > 0$. Assume that $\Phi=\{ \phi_{1}, \ldots 
, \phi_{N} \}$  of ${\cal L}^{\infty}$ is M-refinable and has 
 $L^{2}$-stable shifts. Then we have for any  $1 \leq p, q \leq \infty $, 
 \[
 B^{s,\lambda_{0}}_{pq}(\Phi) \subset B^{s}_{pq}(M)
 \]
provided that there exists a positive number $s_{0}$ with $s_{0} > s$ 
such that $\sup_{l \geq 0}\lambda^{ls_{0}}
|{\rm osc}_{p}^{k}\phi_{j}(\cdot, l)|_{p} <\infty$ 
for all $j=1, \ldots, N$, 
where the norm $|\cdot|_{p}$ and ${\rm osc}_{p}^{k}$ are 
given in {\rm (2)} and {\rm (5)} respectively, and
$\lambda_{0}$ is the least value of absolute values of eigenvalues of $M$.
\end{Proposition}

{\bf Proof}.  
We shall prove for any $f \in B^{s, \lambda_{0}}_{pq}(\Phi)$,
\begin{eqnarray*}
(\sum^{\infty}_{l=0}(\lambda_{0}^{ls} \tilde{\sigma}^{p}_{l}(f))^{q})^{1/q}
 \leq
C(||f||_{p}+(\sum_{l=0}^{\infty}(\lambda_{0}^{ls}\sigma^{p}_{l}
(f))^{q})^{1/q}) 
\end{eqnarray*}
where $\sigma^{p}_{l}$  is the errors of $L^{p}$-approximation given in (10) 
associated with $\Phi$ and $\tilde{\sigma}^{p}_{l}(f)=
||{\rm osc}_{p}^{k}f(\cdot, l)||_{p}$ . 
Since $\sigma^{p}_{l}(f) \rightarrow 0$ as $l \rightarrow \infty \ \ 
(1 \leq p \leq \infty)$, 
we have an $L^{p}$-convergent series 
\[
f(x) =P_{0}f(x) + \sum_{l=0}^{\infty}R_{l}f(x) \equiv 
\sum_{j=1}^{N}\sum_{l=0}^{\infty}
\sum_{\nu \in {\Bbb Z}^{n}}a_{jl}(\nu)\phi_{j}(M^{l}x-\nu)
\]
where $P_{0}f(x)=\sum_{j=1}^{N}\sum_{\nu \in {\Bbb Z}^{n}}
a_{j0}(\nu)\phi_{j}(x-\nu)$ 
and 
$R_{l}f(x)=\sum_{j=1}^{N}\sum_{\nu \in {\Bbb Z}^{n}}
a_{j(l+1)}(\nu)\phi_{j}(M^{l+1}x-\nu)$ are given in (13).

Then we have 
\begin{eqnarray*}
\lefteqn{\tilde{\sigma}^{p}_{l_{0}}(f) 
= \tilde{\sigma}^{p}_{l_{0}}(P_{0}f +\sum_{l=0}^{\infty}R_{l}f)} \\ 
&\leq& \tilde{\sigma}^{p}_{l_{0}}(P_{0}f) + 
\sum_{l=0}^{\infty}\tilde{\sigma}^{p}_{l_{0}}(R_{l}f) 
 \equiv  
I_{0} + \sum_{l=0}^{\infty}I^{\prime}_{l} .
\end{eqnarray*}

We shall give  an estimate of $I_{0}$.
 By (14) we have 
 \begin{eqnarray*}
 \lefteqn{I_{0} \leq 
C\sum_{j=1}^{N}||\sum_{\nu \in {\Bbb Z}^{n}}|a_{j0}(\nu)|
 {\rm osc}_{p}^{k}\phi_{j}(x-\nu , l_{0})||_{p}} \\
 &\leq& 
 C\sum_{j=1}^{N}||a_{j0}||_{l^{p}}|{\rm osc}_{p}^{k}\phi_{j}(\cdot, l_{0})|_{p}
 \leq C||P_{0}f||_{p}\sup_{j}
 |{\rm osc}_{p}^{k}\phi_{j}(\cdot, l_{0})|_{p}.
 \end{eqnarray*}

If $l < l_{0} $, then we see by (14) that 
\begin{eqnarray*}
\lefteqn{I^{\prime}_{l} \leq 
Cm^{-(l+1)/p}\sum_{j=1}^{N}||\sum_{\nu}|a_{j(l+1)}(\nu)|
{\rm osc}_{p}^{k}\phi_{j}(x-\nu, l_{0}-l-1)||_{p}}  \\
&\leq&
C \sum_{j=1}^{N}
 m^{-(l+1)/p}||a_{j(l+1)}||_{l^{p}}|{\rm osc}_{p}^{k}\phi_{j}
 (\cdot, l_{0}-l-1)|_{p}    \leq 
 C||R_{l}f||_{p}\sup_{j}|{\rm osc}_{p}^{k}\phi_{j}(\cdot,l_{0}-l-1)|_{p}.
 \end{eqnarray*}
 
If $l\geq l_{0}$, then we have by the definition,
\[
I^{\prime}_{l} \leq ||R_{l}f||_{p}. 
\]

Hence the above estimates of $I_{0}$ and $I^{\prime}_{l}$, 
 imply that 
\begin{eqnarray*}
\lefteqn{\tilde{\sigma}^{p}_{l_{0}}(f) 
\leq 
C||P_{0}f||_{p}\sup_{j}
|{\rm osc}_{p}^{k}\phi_{j}(\cdot, l_{0})|_{p}} \\
&+&
C\sum_{l=0}^{l_{0}-1}||R_{l}f||_{p}\sup_{j}|{\rm osc}_{p}^{k}\phi_{j}(\cdot, 
l_{0}-l-1)|_{p}+ \sum_{l=l_{0}}^{\infty}||R_{l}f||_{p}
\end{eqnarray*}
and from Hardy's inequality and Theorem 2,
\begin{eqnarray*}
\lefteqn{(\sum_{l_{0}=0}^{\infty}(\lambda_{0}^{l_{0}s}
\tilde{\sigma}^{p}_{l_{0}}(f))^{q})^{1/q}} 
\\ 
&\leq&
C(\sum_{l_{0}=0}^{\infty}(||P_{0}f||_{p}\lambda_{0}^{l_{0}s}
\sup_{j}|{\rm osc}_{p}^{k}\phi_{j}(\cdot, l_{0})|_{p})^{q})^{1/q} \\ 
&+& 
C(\sum_{l_{0}=0}^{\infty}
(\sum_{l=0}^{l_{0}-1}||R_{l}f||_{p}\lambda_{0}^{l_{0}s}
\sup_{j}|{\rm osc}_{p}^{k}\phi_{j}(\cdot, l_{0}-l-1)|_{p})^{q})^{1/q} + 
(\sum_{l_{0}=0}^{\infty}(\sum_{l=l_{0}}^{\infty} 
||R_{l}f||_{p}\lambda_{0}^{l_{0}s})^{q})^{1/q}  \\
&\leq&
 C\sup_{j,l\geq 0}\lambda_{0}^{ls_{0}}|{\rm osc}_{p}^{k}\phi_{j}(\cdot, l)|_{p}
 \{ ||P_{0}f||_{p}( \sum_{l_{0}=0}^{\infty}
 \lambda_{0}^{-l_{0}(s_{0}-s)q})^{1/q} \\ 
 &+& 
(\sum_{l_{0}=0}^{\infty}(\sum_{l=0}^{l_{0}-1}||R_{l}f||_{p}
\lambda_{0}^{l_{0}s}\lambda_{0}^{-(l_{0}-l-1)s_{0}})^{q})^{1/q} \} + 
(\sum_{l_{0}=0}^{\infty}(\sum_{l=l_{0}}^{\infty} 
||R_{l}f||_{p}\lambda_{0}^{l_{0}s})^{q})^{1/q}  \\
&\leq&
C \sup_{j,l\geq 0}\lambda_{0}^{ls_{0}}|{\rm osc}_{p}^{k}\phi_{j}(\cdot, l)|_{p}
 ( ||P_{0}f||_{p} +
(\sum_{l=0}^{\infty}(||R_{l}f||_{p}\lambda_{0}^{ls})^{q})^{1/q}
 ) 
 \\ 
& \leq&
 C\sup_{j,l \geq 0}\lambda_{0}^{ls_{0}}|{\rm osc}_{p}^{k}\phi_{j}
 (\cdot, l)|_{p}(||f||_{p}+(\sum_{l=0}^{\infty}(\sigma_{l}^{p}(f) 
 \lambda_{0}^{ls})^{q})^{1/q}).
 \end{eqnarray*}
 
 This completes the proof of Proposition 2.

A following corollary can be proved by the same way 
in the proof of Proposition 2.

\bigskip

\noindent
{\bf Corollary.}$\ $  {\it 
Given $\lambda >1$ and $s > 0$. Assume that $\Phi=\{ \phi_{1}, \ldots 
, \phi_{N} \}$ and $\Phi^{\prime}=\{ \phi^{\prime}_{1}, \ldots ,
 \phi^{\prime}_{L} \}$ of ${\cal L}^{\infty}$ are M-refinable and have 
 $L^{2}$-stable shifts. Then we have for any  $1 \leq p, q \leq \infty $, 
 \[
 B^{s,\lambda}_{pq}(\Phi^{\prime}) \subset B^{s, \lambda}_{pq}(\Phi)
 \]
provided that there exists a positive number $s_{0}$ with $s_{0} > s$ 
such that $\sup_{l \geq 0}\lambda^{ls_{0}}
|\phi^{\prime}_{j}-P_{l}\phi^{\prime}_{j}|_{p} <\infty$ 
for all $j=1, \ldots, L$, 
where the operator $P_{l}$ is given in {\rm (9)} associated with $\Phi$.
}

\bigskip

For a positive integer $k$ and $1 \leq p \leq \infty$, ${\cal L}_{k}^{p} =
{\cal L}_{k}^{p}({\Bbb R}^{n})$ is denoted to be the space of all functions 
$f$ such that $f(x)(1+|x|)^{k} \in {\cal L}^{p}$.
If $\phi \in L^{p}({\Bbb R}^{n})\ \ (1 \leq p \leq \infty)$ 
is compactly supported, then $\phi \in {\cal L}^{p}_{k}$. 
Furthermore, we observe that if there are 
constants $C>0$ and $\delta >k$ such that 
$|\phi (x)| \leq C(1+|x|)^{-n-\delta}$ 
for all $x \in {\Bbb R}^{n}$
then $\phi \in {\cal L}^{\infty}_{k}$.

For a finite subset $\Phi$ of ${\cal L}^{\infty}_{k}$, the domain of the 
operator $P_{l}$ given in (9), can be extended to include the linear space 
${\Bbb P}^{k}$ of all polynomials of degree no greater 
than $k$ on ${\Bbb R}^{n}$.
For a finite subset $\Phi$ of ${\cal L}^{1}_{k}$, we say that $\Phi$ 
satisfies the Strang-Fix condition of order $k$ if there is a finite 
linear combination $\phi$ of the functions of $\Phi$ and their shifts 
such that 
$\hat{\phi}(0) \neq 0$ and $\partial^{\alpha}\hat{\phi}(2\pi\nu)=0,\ 
|\alpha|\leq k-1,\ \nu \in {\Bbb Z}^{n}$ with $\nu \neq 0$.

\bigskip

\noindent
\begin{Lemma}.  
Let $\Phi$ be a finite subset of ${\cal L}^{\infty}_{k}$ that has $L^{2}$-
stable shifts. Then 
 $\Phi$ satisfies the Strang-Fix condition of order $k$ if and only if 
 $P_{0}q = q$ for any $q \in {\Bbb P}^{k-1}$.

Moreover, if this is the case, then we have 
 $||P_{l}f -f ||_{p} \leq C \lambda_{0}^{-lk}\sum_{|\alpha|=k}
||\partial^{\alpha}f||_{p}$ for any $f$ in the Sobolev space 
$W^{p}_{k}({\Bbb R}^{n})\ \ (1 \leq p \leq \infty)$, 
with a constant $C$ independent of $f,p$ and $l$ 
 where $\lambda_{0}$ is the least value of absolute values of eigenvalues 
 of the dilation matrix $M$, that is, $W^{p}_{k}({\Bbb R}^{n}) \subset 
 B^{s, \lambda_{0}}_{pq}(\Phi)$ if $0<s<k$ and $1 \leq q \leq \infty$.
\end{Lemma}

{\bf Proof}. We can prove by the same way of [8,  Theorem 5.2]. 
We will omit its details.

\section{Characterization of Besov spaces}

\bigskip

Let $\Pi$ be a self-affine lattice tiling 
$\{T+\nu \}_{\nu \in {\Bbb Z}^{n}}$ and $\Pi_{l}$ denote the subdivision 
$\{ M^{-l}(T+\nu) \}_{\nu \in {\Bbb Z}^{n}}$ of ${\Bbb R}^{n}$ for 
a nonnegative integer $l$. Let $\Phi=\{ \phi_{1}, \ldots, 
\phi_{N} \}$ be a finite subset of ${\cal L}^{\infty}$ and $\lambda_{0}$ the 
least value of absolute values of eigenvalues of the dilation matrix $M$.

\begin{Proposition}. 
Given $1 \leq p,\ q \leq \infty$ and $k>s>0$. 
Assume that a finite subset $\Phi=\{ \phi_{1}, \ldots, \phi_{N} \}$ 
of ${\cal L}^{\infty}_{k}$ satisfies

{\rm (a)}   
 $\Phi$ has $L^{2}$-stable shifts,

{\rm (b)}  $\Phi$ is $M$-refinable,

{\rm (c)} $\Phi$ satisfies the Strang-Fix condition of order $k$. 

Then we have 
$
B^{s}_{pq}(M) \subset B^{s, \lambda_{0}}_{pq}(\Phi)$.

\end{Proposition}

{\bf Proof}.  
 We shall prove for any $f \in B^{s}_{pq}(M)$ by the same routine of the 
 proof of Theorem 1,
 \begin{eqnarray*}
 (\sum_{l=0}^{\infty}(\lambda_{0}^{ls}\sigma_{l}^{p}(f))^{q})^{1/q} 
 \leq C||f||_{B^{s}_{pq}(M)} 
 \end{eqnarray*}
 where $\sigma^{p}_{l}$ is given in (10) associated with $\Phi$.
 We use the same notations in Theorem 1. 
  We choose a function $\chi$ in $C^{\infty}_{c}({\Bbb R}^{n})$ such that 
 $\int |\chi(u)| du =1$ and supp $\chi \subset \{ u \in {\Bbb R}^{n}:  
 |u| < r/2k \}$ where $r$ is the positive number given in (7). We write 
 $\chi_{l}(u)=m^{l}\chi(M^{l}u),\ \ h_{l}(x)=\int(f(x)-\triangle_{u}^{k}f(x))
 \chi_{l}(u)du$ and $g_{l}=P_{l}h_{l}-h_{l}$
  where $P_{l}$ is given in (9) associated with $\Phi$. 
  Then we have for $1 \leq p \leq 
 \infty$, 
 \[
 ||f-P_{l}f||_{p} \leq ||f-h_{l}||_{p} + ||g_{l}||_{p} 
 +||P_{l}h_{l}-P_{l}f||_{p} \leq C||f-h_{l}||_{p} + ||g_{l}||_{p} 
 \equiv CI_{1}+I_{2}.
 \]
 Obviously we have :
 \[
 I_{1} \leq C\sup_{k|M^{l}u|< r/2}||\triangle^{k}_{u}f||_{p}. 
 \]
We shall give an estimate of $I_{2}$ by (1): 
\begin{eqnarray}
I_{2} = (\sum_{Q \in \Pi_{l}}\int_{Q}|g_{l}(x)|^{p}dx)^{1/p}
=(\sum_{\nu \in {\Bbb Z}^{n}}\int_{M^{-l}T}
|g_{l}(x-M^{-l}\nu)|^{p}dx)^{1/p}.
\end{eqnarray}
Let $q_{z}$ be the $(k-1)$-th Taylor polynomial of $h_{l}$ about 
$z \in {\Bbb R}^{n}$ and let $r_{z}$ be the corresponding remainder. 
Since $\Phi$ satisfies the Strang-Fix condition of order $k$, 
we see from Lemma 1 
\[
g_{l}(x-M^{-l}\nu)=P_{l}r_{x-M^{-l}\nu}(x-M^{-l}\nu)=
m^{l}\int K(M^{l}x, M^{l}y)r_{x-M^{-l}\nu}(y-M^{-l}\nu)dy
\]
where $K(x,y)=\sum_{j=1}^{N}\sum_{\nu \in {\Bbb Z}^{n}}
\phi_{j}(x-\nu)\bar{\tilde{\phi}}_{j}(y-\nu)$.

To estimate  $I_{2}$, we use
\[
r_{x-M^{-l}\nu}(y-M^{-l}\nu) = 
\int_{0}^{1}\sum_{|\beta|=k}\displaystyle{\frac{k}{\beta!}}
\partial^{\beta}h_{l}(x+t(y-x)-M^{-l}\nu)(1-t)^{k-1}(y-x)^{\beta}dt,
\]
and 
\begin{eqnarray*}
\lefteqn{|\partial^{\beta}h_{l}(x)|  
\leq 
C\sum_{e=1}^{k}
(\int_{|u|< r/2k}|f(x-eM^{-l}u)|^{p}du)^{1/p}} \\
 &\leq&
  C\sum_{e=1}^{k}(m^{l}\int_{|M^{l}u|< re/2k}|f(x-u)|^{p}du)^{1/p} 
\leq
  Cm^{l/p}(\int_{|M^{l}u|< r/2}|f(x-u)|^{p}du)^{1/p}.
\end{eqnarray*}
Hence we get an estimate:
\begin{eqnarray*}
\lefteqn{
(\sum_{\nu \in {\Bbb Z}^{n}}
|r_{x-M^{-l}\nu}(y-M^{-l}\nu)|^{p})^{1/p}
}
 \\ 
&\leq& 
 C\int_{0}^{1}\sum_{|\beta|=k}(\sum_{\nu}
|\partial^{\beta}h_{l}(x+t(y-x)-M^{-l}\nu)|^{p})^{1/p}
(1-t)^{k-1}|x-y|^{k}dt
 \\
&\leq& 
C\int_{0}^{1}\sum_{|\beta|=k}(\sum_{\nu}m^{l}\int_{|M^{l}u|< r/2}
|f(x+t(y-x)-M^{-l}\nu-u)|^{p}du)^{1/p}(1-t)^{k-1}|x-y|^{k}dt
 \\
&\leq& 
C\int_{0}^{1}m^{l/p}(\sum_{\nu}\int_{M^{-l}(T+\nu)} 
|f(x+t(y-x)+u)|^{p}du)^{1/p}(1-t)^{k-1}|x-y|^{k}dt 
 \\
&\leq& 
C\int_{0}^{1}m^{l/p}||f||_{p}(1-t)^{k-1}|x-y|^{k}dt 
\leq C|x-y|^{k}m^{l/p}||f||_{p}
.
\end{eqnarray*}
Hence, since $\Phi \subset {\cal L}^{\infty}_{k}$,
 we get an estimate of $I_{2}$ in (17): 
 \begin{eqnarray*}
 \lefteqn{I_{2}  
\leq 
Cm^{l}(\int_{M^{-l}T}\sum_{\nu}
(\int|K(M^{l}x, M^{l}y)||r_{x-M^{-l}\nu}(y-M^{-l}\nu)|dy)^{p}dx)^{1/p}} 
\\
&\leq& 
Cm^{l}(\int_{M^{-l}T}
(\int|K(M^{l}x, M^{l}y)|(\sum_{\nu}
|r_{x-M^{-l}\nu}(y-M^{-l}\nu)|^{p})^{1/p}dy)^{p}dx)^{1/p} \\
&\leq& 
 Cm^{l+l/p}||f||_{p}
 (\int_{M^{-l}T}(\int|K(M^{l}x, M^{l}y)||x-y|^{k}dy)^{p}dx)^{1/p} \\
&\leq & 
C||f||_{p}(\int_{T}(\int |K(x, y)||M^{-l}(x-y)|^{k}dy)^{p}dx)^{1/p} 
\\ 
&\leq& 
C||f||_{p}\lambda_{0}^{-lk}(\int_{T}(\int |K(x, y)||x-y|^{k}dy)^{p}dx)^{1/p} 
 \leq 
 C||f||_{p}\lambda_{0}^{-lk}.
\end{eqnarray*}
Now we combine the estimates of $I_{1}$ and $I_{2}$ to write 
\begin{eqnarray*}
||f-P_{l}f||_{p} \leq CI_{1}+I_{2} \leq 
C(\sup_{k|M^{l}u|< r/2}||\triangle^{k}_{u}f||_{p}+
\lambda_{0}^{-lk}||f||_{p}).
\end{eqnarray*}
This implies that
\[
(\sum_{l=0}^{\infty}(\lambda_{0}^{ls}\sigma_{l}^{p}(f))^{q})^{1/q}
\leq C||f||_{B^{s}_{pq}(M)}.
\]
This completes the proof of Proposition 3.

\bigskip

Even if $\Phi$ in Proposition 3 is $M$-refinable for a.e. 
$x \in {\Bbb R}^{n}$, Proposition 3 is true. Hence we have a following 
corollary:

\bigskip

\noindent
{\bf Corollary}. {\it 
Given $1 \leq p, q \leq \infty$ and $0 < s < k$. Then 
\[
B^{s}_{pq}(M) \subset 
B^{s, \lambda_{0}}_{pq}(\{p_{\alpha}\}_{|\alpha|< k})
\]
where the functions $p_{\alpha}$ are given in {\rm (3)}.
}

\bigskip

\noindent
{\bf Remark 3}. 
(a) We can define the operator $P_{\Pi_{l}}(l=0,1,2,\ldots)$ associated with 
$\Phi=\{ p_{\alpha} \}_{|\alpha|\leq k}$, given in (9). 
Then we have $P_{\Pi_{l}}f(x)=\sum_{Q \in \Pi_{l}}P_{Q}f(x)$ where 
\[P_{Q}f=\sum_{|\alpha| \leq k}\langle f, {\tilde p}_{\alpha}
(M^{l}\cdot-x_{0})\rangle m^{l}p_{\alpha}(M^{l}x-x_{0})
\]  
for $Q=M^{-l}(T+x_{0})$ is of type . 
We denote by $\tilde{\omega}_{p}^{k}(f,\Pi_{l})$ the error  of 
$L^{p}$-approximation in (10) associated with 
$\Phi=\{ p_{\alpha}\}_{|\alpha|\leq k}$.
We can prove a following equivalences for $f \in L^{p}({\Bbb R}^{n})
\ (1 \leq p \leq \infty)$ by using the results in [2] and [7]:
\[
\tilde{\omega}^{k}_{p}(f,\Pi_{l})
\sim (\sum_{Q \in \Pi_{l}}\int_{Q}|f-P_{Q}f|^{p}dy)^{1/p} 
\sim 
(\sum_{Q \in \Pi_{l}}\inf_{P \in {\Bbb P}^{k}}\int_{Q}|f-P|^{p}dy)^{1/p} .
\]

\noindent
(b)  Let $\Pi_{l}(x_{0})$ denote $\{ M^{-l}(T+\gamma)+x_{0} \}
_{\gamma \in {\Bbb Z}^{n}}$ for $x_{0} \in {\Bbb R}^{n}$. 
We write 
\[
\tilde{\omega}^{k}_{p}(f, \Pi_{l}(x_{0}))=(\sum_{Q \in \Pi_{l}(x_{0})}
\int_{Q}|f(y)-P_{Q}f(y)|^{p}dy)^{1/p}
\]
 and 
$\tilde{\omega}^{k}_{p}(f,l)=\sup _{x_{0} \in {\Bbb R}^{n}}
\tilde{\omega}^{k}_{p}(f , \Pi_{l}(x_{0}))$.
When the self-affine lattice tiling $\Pi$ is the net of closed cubes 
generated by $T=[0,1]^{n}$ and the dilation matrix $M$ is $2Id$, 
we see
for $1 \leq p,q \leq \infty$ and $k+1 >s > 0$
\[
||f||_{B^{s}_{pq}}({\Bbb R}^{n}) \sim 
||f||_{p} + (\sum_{l=0}^{\infty}(\lambda_{0}^{ls}\tilde{\omega}^{k}_{p}(f, l))
^{q})^{1/q}.
\]
(See [7]).

\bigskip

\begin{Theorem}.
Given  $1 \leq p,\ q \leq \infty$ and $k>s>0$. 
Assume that a finite subset $\Phi=\{ \phi_{1}, \ldots, \phi_{N} \}$ 
of ${\cal L}^{\infty}_{k}$ satisfies

{\rm (a)}   
 $\Phi$ has $L^{2}$-stable shifts,

{\rm (b)}  $\Phi$ is $M$-refinable,

{\rm (c)} there exists a positive number $s_{0}$ with $s_{0}>s$ 
such that $\ \sup_{l\geq 0}\lambda_{0}^{ls_{0}}
|{\rm osc}_{p}^{k-1}\phi_{j}(\cdot, l)|_{p}<\infty$ for all $j=1, \ldots,N,$

{\rm (d)} $\Phi$ satisfies the Strang-Fix condition of order $k$. 

Then we have  
$B^{s}_{pq}(M)=B^{s, \lambda_{0}}_{pq}(\Phi)$  with 
 equivalent norms   
\[
||f||_{B^{s}_{pq}(M)} 
\sim 
||f||_{B^{s,\lambda_{0}}_{pq}(\Phi)} 
\]
where the norms $||f||_{B^{s}_{pq}(M)}$ and 
$||f||_{B^{s, \lambda_{0}}_{pq}(\Phi)}$
 are given in {\rm (6)} and {\rm (11)} respectively, 
 and $\lambda_{0}$ is the least value of absolute values of eigenvalues of 
 the dilation matrix $M$.
\end{Theorem}

\bigskip

\noindent
{\rm\bf Remark 4}.   When $\{ \phi_{j} \}_{j=1}^{N}$ have compact supports, we 
see that the condition (c) in Theorem 3 can be rephrased as :

${\rm (c)^{\prime}}$ There exists a positive number $s_{0}>s$ such that 
$\sup_{l \geq 0}\lambda_{0}^{ls_{0}}||{\rm osc}_{p}^{k-1}
\phi_{j}(\cdot, l)||_{p}
< \infty$, (that is, $\phi_{j} \in B^{s_{0}}_{p\infty}(M)$ 
if $s_{0} < k$) for all $j=1, \ldots, N$.

\bigskip

\noindent
{\rm\bf Proof of Theorem 3.} 
This result is an immediate consequence of Proposition 2 
and Proposition 3.

\bigskip

We say that a function on ${\Bbb R}^{n}$ is $k$-regular 
if it is of class $C^{k}$ and 
rapidly decreasing in the sense that 
$|\partial^{\alpha}f(x)| \leq C_{N}(1+|x|)^{-N}$  
for all $N=0,\ 1, \ 2,\ \ldots$ and all $|\alpha| \leq k$. 
Any $k$-regular function belongs to ${\cal L}^{\infty}_{N}$ 
for any $N \geq 0$ and any $k$-regular function $f$ satisfies the condition 
(c) in Theorem 3 : $\sup_{l \geq 0}\lambda_{0}^{lk}|{\rm osc}_{p}^{k-1}
f(\cdot, l)|_{p} < \infty$.

\bigskip

\begin{Corollary}.
Suppose that a dilation matrix is of the form $M=\lambda_{0}Id$ 
with $\lambda_{0} >1$. 
Let $1 \leq p,\ q \leq \infty$ and $k>s>0$.  
Assume that a finite subset $\Phi=\{ \phi_{1}, \ldots , \phi_{N} \}$ of 
$k$-regular functions on 
${\Bbb R}^{n}$ satisfies:

{\rm (a)}   
 $\Phi$ has $L^{2}$-stable shifts,

{\rm (b)}  $\Phi$ is $M$-refinable.

Then there exits a set $\{ \psi^{\epsilon}_{j} \}_{j=1, \ldots, N}^{\epsilon=1, \ldots, m-1}$ of $k$-regular wavelets  associated with  $\Phi$, and 
 we have equivalence of norms, if one of them exit,
for wavelet expansion given in {\rm (15)} 
with the usual modification for $q=\infty$: 
\begin{eqnarray*}
&{\rm (i)}& \lefteqn{ \ ||f||_{B^{s}_{pq}(M)},} \\ 
\sim 
&{\rm (ii)}& 
||f||_{B^{s,\lambda_{0}}_{pq}(\Phi)}, \\
\sim 
&{\rm (iii)}&  
 \sum_{j=1}^{N}||a_{j0}||_{l^{p}} + 
(\sum_{l=0}^{\infty}(\lambda_{0}^{l(s+n/2-n/p)}
\sum_{j=1}^{N}\sum_{\epsilon=1}^{m-1}||b_{jl}^{\epsilon}||_{l^{p}})^{q})^{1/q}. \end{eqnarray*}
\end{Corollary}

\bigskip

{\rm\bf Proof}.$\ \ $  
 From [15,  Theorem 5.15], for a finite subset $\Phi$ of 
 $k$- regular functions
 there exists an associated set of $k$-regular wavelets 
 for a general dilation matrix $M$ if $m > (n+1)/2$. 
Since a finite subset of $k$-regular functions satisfies 
the Strang-Fix condition of order $k+1$ in the case $M=\lambda_{0}Id$  
(See  [9,  Theorem 4 in 2.6] and Lemma 1), 
we have  the equivalence of (i) and (ii) from Theorem 3. 
The equivalence of (ii) and (iii)  can be proved by Theorem 2.

\bigskip

We define the tensor product B-spline by 
${\cal M}_{k}=\prod_{i=1}^{n}M_{k}(x_{i}), \ \ x=(x_{1}, \ldots, x_{n})
 \in {\Bbb R}^{n},\ \ k=1,2, \ldots.$ 
 where $M_{k}(t)$ is the $k$-th order central B-spline, that is, 
 $\hat{M}_{k}(t)=(\displaystyle{\frac{\sin (t/2)}{t/2}})^{k}$. 
 Let us denote by $\{ e^{i} \}_{i=1}^{n}$ the set of unit vectors 
 in ${\Bbb R}^{n}$. We put $e^{n+1}=\sum_{i=1}^{n}e^{i}$, and 
 $X=\{ x^{1}, \ldots, x^{d_{0}} \}$ with $x^{1}=e^{1}, \ldots, x^{d_{1}}=e^{1},
  x^{d_{1}+1}=e^{2}, \ldots, x^{d_{1}+d_{2}}=e^{2}, \ldots, 
 x^{d_{1}+\cdots+d_{n}+1}=e^{n+1}, \ldots, x^{d_{0}}=e^{n+1}$ where 
 $d_{0}=d_{1}+ \cdots +d_{n+1}$. 
 We denote the box spline $B(x,X)$ corresponding to $X$ given by 
 $\hat{B}(x,X)=(2\pi)^{-n/2}\prod_{j=1}^{d_{0}}\displaystyle
 {\frac{1-e^{ix^{j}\cdot x}}{ix^{j}\cdot x}}$.
  In the case that the self-affine lattice tiling is the net of closed cubes 
 generated by  $T=[0, 1]^{n}$  and the dilation matrix is $2Id$, 
   the $k$-th order tensor product 
  B-spline ${\cal M}_{k}$ satisfies  the conditions of  Theorem 3, 
  particularly, 
  ${\cal M}_{k} \in B^{k-1+1/p}_{p\infty}({\Bbb R}^{n})$   and 
  ${\cal M}_{k}$ satisfies the Strang-Fix condition of order $k$. 
  The above box spline $B(x , X)$ also satisfies the conditions of 
  Theorem 3 replacing the above $k$ by $k=\min\{ d_{i}+d_{j}
  : \ \ i,j =1,\ldots, n+1,\ \ i \neq j \}$. Hence we get results 
  of [3] and [12].

\begin{Corollary}. 
Suppose that the self-affine lattice tiling is the net 
$\Pi=\{ T+\nu \}_{\nu \in {\Bbb Z}^{n}}$ of closed cubes 
generated by $T=[0,1]^{n}$ 
 and the dilation matrix is $2Id$. Then Theorem 3 remains true 
for the tensor product B-spline $\Phi=\{ {\cal M}_{k} \}$ or 
the box spline $\Phi=\{ B(x,X) \}$.

\end{Corollary}

More general results have been given in W.Sickel [12] and R.A. DeVore , 
B. Jarwerth and V. Popov [4] 
in a case of some compactly supported functions. 
 Theorem 3 is a generalization of those results.

A following proposition is a pointwise version 
of Corollary 1 in Theorem 3.

\bigskip

\begin{Proposition}.
Suppose that a dilation matrix is of the form $M=\lambda_{0}Id$ 
with $\lambda_{0} >1$ and $k > s >0$.  
Assume that a finite subset $\Phi=\{ \phi_{1}, \ldots , \phi_{N} \}$ of 
$k$-regular functions on 
${\Bbb R}^{n}$ satisfies:

{\rm (a)}   
 $\Phi$ has $L^{2}$-stable shifts,

{\rm (b)}  $\Phi$ is $M$-refinable.

Then for $x \in {\Bbb R}^{n}$ and a bounded function f on ${\Bbb R}^{n}$ ,
 following properties are equivalent: 
\begin{eqnarray*}
&{\rm (i)}&  f \in C^{s}(x) , \\
&{\rm (ii)}& 
|f(y)-P_{l}f(y)| \leq C(\lambda_{0}^{-l}+|x-y|)^{s}\ \ \ l \geq 0 
\end{eqnarray*}
where $P_{l}f$ is given in {\rm (9)}.

 \end{Proposition}

{\bf Proof}. This can be proved by the same way as in Proposition 1. 
 See [1, Theorem 3].

\bigskip

\noindent
{\bf Corollary}. {\it 
Suppose that the conditions in Proposition 4 are satisfied.  
Let $s>s^{\prime}>0$.

{\rm (a)} 
If $f \in C^{s}(x)$, we have 
\[
|R_{l}f(y)| \leq C(\lambda_{0}^{-l}+|x-y|)^{s} \ \ \ l=0,1,2,3,\ldots
\]
where $R_{l}f$ is given in {\rm (12)}.

If it holds  
\[
|R_{l}f(y)| \leq C\lambda_{0}^{-sl}(1+\lambda_{0}^{l}|x-y|)^{s^{\prime}}
\ \ \ l=0, 1, 2, 3, \ldots ,   
\]
then $f \in C^{s}(x)$.

{\rm (b)} 
If $f \in C^{s}(x)$, we have 
\[
|b^{\epsilon}_{jl}(\nu)| \leq C\lambda_{0}^{-(s+\frac{n}{2})l}
(1+|\lambda_{0}^{l}x-\nu|)^{s} 
\]
for $j=1,\ldots, N, l=1, 2, 3,\ldots, \epsilon =1 , \ldots, m-1$ 
and any $\nu \in {\Bbb Z}^{n}$ where $b^{\epsilon}_{jl}(\nu)$ 
is given in {\rm (16)}.

If it holds  
\[
|b^{\epsilon}_{jl}(\nu)| \leq C\lambda_{0}^{-(s+\frac{n}{2})l}
(1+|\lambda_{0}^{l}x-\nu|)^{s^{\prime}}
\ \ {\rm for} \ \ j=1, \ldots, N,\ \ 
l=1, 2, 3, \ldots \ {\rm and} \ \ \epsilon =1, \ldots, m-1 
\]
and any $ \nu \in {\Bbb Z}^{n}$,    
then $f \in C^{s}(x)$.

{\rm (c)} For $\{ a_{jl}(\nu) \}$ given in {\rm (13)}, 
 if it holds 
\[
|a_{jl}(\nu)| \leq C\lambda_{0}^{-sl}(1+|\lambda_{0}^{l}x-\nu|)^{s^{\prime}} 
\ \ j=1, \ldots, N, \ \ l > 0  \ {\rm and}\ \ \nu \in {\Bbb Z}^{n},
\]
then $f \in C^{s}(x)$.
}

\section{Scaling exponents}
 For $1 \leq p,q \leq \infty$ we define  
$\alpha_{pq}(f)=\sup \{ s \geq 0 : f \in B^{s}_{pq}(M) \}$ 
for functions $f \in L^{p}({\Bbb R}^{n})$. 
If there is not a positive number $s$ with $f \in B^{s}_{pq}(M)$, then we 
define $\alpha_{pq}(f)=0$. We remark that $\alpha_{pq}(f) > 0$ 
for any $f \in L^{p}({\Bbb R}^{n})$ in the case $1 \leq p <\infty$. In the 
same manner we define $\alpha_{pq}(f,x) = 
\sup \{s \geq 0: f \in T^{s}_{pq}(x) \}$ 
for $x \in {\Bbb R}^{n}$ and bounded functions $f$ on ${\Bbb R}^{n}$. 
We put $ \alpha_{p}(f)=\alpha_{p \infty}(f)$, $\alpha(f)=\alpha_{\infty}(f)$, 
$\alpha_{p}(f,x)=\alpha_{p\infty}(f,x)$ and $\alpha(f,x)=\alpha_{\infty}(f,x)$.

We can prove a following proposition by the embedding theorem ( See 
Remark 1, 2 and [11]).
\bigskip

\begin {Proposition}.

{\rm (i)} $\alpha_{p}(f)=\alpha_{p\eta}(f)$ 
for $1 \leq p, \eta \leq \infty$,

{\rm (ii)} $\alpha(f) > \alpha_{p}(f) -\frac{n}{p} \geq 
\alpha_{q}(f) - \frac{n}{q}$ for $1 \leq q \leq p < \infty$ 
when $M=\lambda_{0}Id$,

{\rm (iii)} $\alpha_{p}(f, x) = \alpha_{p\eta}(f,x)$ 
for $1 \leq p, \eta \leq \infty$,

{\rm (iv)} $\alpha(f) \leq \alpha(f,x) \leq \alpha_{p}(f,x) 
\leq \alpha_{q}(f,x)$ for $1 \leq q \leq p < \infty.$
\end{Proposition}

\bigskip

For $1 \leq p \leq \infty$ we have by Theorem 1 
and Theorem B 
\[
\alpha_{p}(f)=-\displaystyle{\frac{\log A_{p}(f)}{\log \lambda_{0}}}
\]
if the right hand side of the above equality is less than $k+1$ where 
\[
A_{p}(f)
=\limsup_{l \rightarrow \infty} ||{\rm osc}^{k}_{p}f(\cdot, l)||_{p}^{1/l}
=\limsup_{l \rightarrow \infty} 
\sup_{(k+1)|M^{l}u|<r/2}||\triangle^{k+1}_{u}f||_{p}^{1/l}
\]
 and furthermore when $M=\lambda_{0}Id$ with $\lambda_{0} > 1$
\[
A_{p}(f)
=\limsup_{l \rightarrow \infty} ||f-S_{l}f||_{p}^{1/l}
=\limsup_{l \rightarrow \infty} ||f_{l}||_{p}^{1/l}.
\]
For $1 \leq p \leq \infty$ we have by the corollary 
of Theorem 1  
\[
\alpha_{p}(f,x)=-\displaystyle{\frac{\log A_{p}(f,x)}{\log \lambda_{0}}}
\]
if the right hand side of the above equality is less than $k+1$ where 
\[
A_{p}(f,x)
=\limsup_{l \rightarrow \infty} {\rm osc}^{k}_{p}f(x , l)^{1/l}
=\limsup_{l \rightarrow \infty} 
\sup_{(k+1)|M^{l}u|<r/2}(\displaystyle{\frac{1}{|Q_{l}(x)|}}\int_{Q_{l}(x)}
|\triangle^{k+1}_{u}f(y)|^{p}dy)^{1/pl}.
\]
Furthermore when $M=\lambda_{0}Id$ with $\lambda_{0} > 1$, 
we have by Proposition 1 and its corollary 
\[
\alpha(f,x)
 =  \liminf_{\lambda_{0}^{-l}+|x-y| \rightarrow 0}
\displaystyle{\frac{\log |f(y)-S_{l}f(y)|}
{\log (\lambda_{0}^{-l}+|x-y|)}} 
\]
and, if $\alpha(f) > 0$ 
\[
\alpha(f,x)  =  \liminf_{\lambda_{0}^{-l}+|x-y| \rightarrow 0}
\displaystyle{\frac{\log |f_{l}(y)|}
{\log (\lambda_{0}^{-l}+|x-y|)}} 
\]
where $S_{l}f$ and $f_{l}$ are given 
for Littlewood-Paley decompostion in (8).

We can prove a following proposition by Theorem 2, Theorem 3, Proposition 4 
and its corollary.

\bigskip

\begin{Proposition}. 
{\rm (i)}.  
Assume that a finite subset $\Phi=\{  \phi_{1}, \ldots , \phi_{N} \}$ 
of ${\cal L}^{\infty}_{k}$ satisfies 
the conditions {\rm (a), (b), (c)} and {\rm (d)} of  Theorem 3.

 Then for $f \in L^{p}({\Bbb R}^{n})\ 
(1 \leq p \leq \infty)$ we have 
\[
\alpha_{p}(f)=-\displaystyle{\frac{\log A_{p}(f)}{\log \lambda_{0}}}
= \displaystyle{\frac{\log m}{p\log \lambda_{0}}}-
\displaystyle{\frac{\log \rho_{p}(f)}{\log \lambda_{0}}}
\]
if the second and third parts of the above equality are less than 
{\rm min(}k, $s_{0}${\rm )} where 
\[
A_{p}(f)
= \limsup_{l \rightarrow \infty} \sigma^{p}_{l}(f)^{1/l}
=\limsup_{l \rightarrow \infty} ||R_{l}(f)||_{p}^{1/l}
\] and    
\[
\rho_{p}(f)=\limsup_{l \rightarrow \infty}\sum_{j=1}^{N}||a_{jl}||
_{l^{p}}^{1/l} 
=\inf \limsup_{l \rightarrow \infty} \sum_{j=1}^{N}
||c_{jl}||_{l^{p}}^{1/l}
\]
and $\{ a_{jl} \}$ is given by {\rm (13)} 
and {\rm inf} is taken over all admissible representations 
$f(x) = \sum_{j=1}^{N}\sum_{l=0}^{\infty}\sum_{\nu \in {\Bbb Z}^{n}}
c_{jl}(\nu)\phi_{j}(M^{l}x-\nu)$ as in Theorem 2.

{\rm (ii)}.  Furthermore  when  $m > (n+1)/2$,
 we have 
\[
\alpha_{p}(f)=(1/p-1/2)\displaystyle{\frac{\log m}{\log \lambda_{0}}}-
\displaystyle{\frac{\log \rho^{\prime}_{p}(f)}{\log \lambda_{0}}}
\]
if the right hand side of the above equality is less than 
 {\rm min(}k, $s_{0}${\rm )} where 
\[
\rho^{\prime}_{p}(f) =\limsup_{l \rightarrow \infty}\sum_{j=1}^{N}
\sum_{\epsilon=1}^{m-1}||b_{jl}^{\epsilon}||_{l^{p}}^{1/l}
\]
and  $\{ b_{jl}^{\epsilon} \}$ is given in {\rm (16)} 
for the wavelet expansion {\rm (15)} associated with $\Phi$.

{\rm (iii)}. 
Suppose that conditions in Proposition 4 hold for a bounded function $f$. 
Then  we have 
\[
\alpha(f,x)
 =  \liminf_{\lambda_{0}^{-l}+|x-y| \rightarrow 0}
\displaystyle{\frac{\log |f(y)-P_{l}f(y)|}
{\log (\lambda_{0}^{-l}+|x-y|)}} 
\]
if the right hand side of the above equality is less than $k$ 
and, 
\begin{eqnarray*}
\lefteqn{
\alpha(f,x) =  \liminf_{\lambda_{0}^{-l}+|x-y| \rightarrow 0}
\displaystyle{\frac{\log |R_{l}f(y)|}
{\log (\lambda_{0}^{-l}+|x-y|)}} 
} \\ 
&=&
\liminf_{\lambda_{0}^{-l}+|x-\lambda_{0}^{-l}\nu| \rightarrow 0}\inf_{j}
\displaystyle{\frac{\log \lambda_{0}^{\frac{n}{2}l}|b_{jl}^{\epsilon}(\nu)|}
{\log (\lambda_{0}^{-l}+|x-\lambda_{0}^{-l}\nu|)}} 
\leq  \liminf_{ \lambda_{0}^{-l}+|x-\lambda_{0}^{-l}\nu|\rightarrow 0}
\inf_{j}\displaystyle{\frac{\log |a_{jl}(\nu)|}
{\log (\lambda_{0}^{-l}+|x-\lambda_{0}^{-l}\nu|)}} 
\end{eqnarray*}
if $\alpha(f)>0$ and the right hand side of the above inequality is less 
than $k$ where $P_{l}f$, $R_{l}f$  and $\{ a_{jl} \}$ are given 
in {\rm (9)}, {\rm (12)} and {\rm (13)} respectively.  
\end{Proposition}

\bigskip

Let $\Pi=\{ T+\nu \}_{\nu \in {\Bbb Z}^{n}}$ be a self-affine lattice tiling 
with a dilation matrix $M$ and a set $\Gamma_{0}$ of digits, 
and $\Pi_{l}$ denote the subdivision 
$\{  M^{-l}(T+\nu) \}_{\nu \in {\Bbb Z}^{n}}$ of ${\Bbb R}^{n}$ 
for a nonnegative integer $l$. We write  $Q= M^{-l}(T+\nu_{Q})$ for 
$Q \in \Pi_{l}$. Let $\Pi_{l}(T)=\{ Q \in \Pi_{l}: \ Q \subset T \}$ and 
$\Pi(T)=\cup_{l=0}^{\infty}\Pi_{l}(T)$. We put 
$\Gamma_{0}=\{ \gamma_{1}, \cdots , \gamma_{m} \}$. Then from (1) for 
$Q \in \Pi_{l}(T)$, $\nu_{Q}$ is of a form 
$\nu_{Q}= M^{l-1}\gamma_{i_{1}} + \cdots +\gamma_{i_{l}}, \ \ 
\gamma_{i_{1}}, \cdots , \gamma_{i_{l}} \in \Gamma_{0}$ and let 
$M_{Q}y=M^{l}y-\nu_{Q}$ and $\mu_{Q} = \mu_{i_{1}}\cdots \mu_{i_{l}}$ 
for $l > 0$ where $\mu_{1}, \mu_{2}, \ldots , \mu_{m}$ are real 
or complex numbers with $0 < |\mu_{i}|< 1$ , $i=1, \ldots, m$. For $l=0$ 
we put $M_{T}=Id$ and $\mu_{T}=1$.

From now we suppose that a dilation matrix $M$ 
is of a form $M=\lambda_{0} Id$ with 
$\lambda_{0}>1$ and we consider a bounded function $f$ 
which is given by a series 
\begin{eqnarray}
f(y)=\sum_{Q \in \Pi(T)} \mu_{Q} \phi(M_{Q}y), \ \ y \in {\Bbb R}^{n}
\end{eqnarray}
where a function $\phi$ is bounded and zero outside $T^{o}$. 
We remark that  $\alpha(f) \leq \alpha(\phi)$. 
Let 
\[
\tau_{0}(x)\equiv 
\liminf_{l \rightarrow \infty} \inf_{K_{l}(x) \ni Q} 
\displaystyle{\frac{\log |\mu_{Q}|}
{\log (\lambda_{0}^{-l}+|x-\lambda_{0}^{-l}\nu_{Q}|)}}
=\liminf_{l \rightarrow \infty} \inf_{K_{l}(x) \ni Q} 
\displaystyle{\frac{\log |\mu_{Q}|}
{\log \lambda_{0}^{-l}}}
\]
where $K_{l}(x) \equiv \{ Q \in \Pi_{l}(T):\  B(x, \lambda_{0}^{-l}) 
\cap Q \neq \emptyset \}$ and $B(x, \lambda_{0}^{-l})$ is a ball centered 
at $x$ with a radius $\lambda_{0}^{-l}$. When $x \in \Omega \equiv 
\cap_{l=0}^{\infty} \cup_{Q \in \Pi_{l}(T)}Q^{o}$(the interior of $Q$) 
there exits a unique sequence $\{ Q_{l,x} \}_{l \geq 0}$ such that 
$Q_{l,x} \in \Pi_{l}(T)$ and $x \in Q_{l,x}^{o}$. 
Then we have for $x \in \Omega$
\[
\tau_{0}(x)=\liminf_{l \rightarrow \infty}\displaystyle{\frac
{\log |\mu_{Q_{l,x}}|}
{\log \lambda_{0}^{-l}}}.
\]
Let for $x \in \Omega$ 
\[
\tau_{1}(x)\equiv\liminf_{l \rightarrow \infty}
\displaystyle{\frac{\log |\mu_{Q_{l,x}}|}{\log \Delta_{l}(x)}}
\]
where $\Delta_{l}(x)= {\rm dist}(x, \partial Q_{l,x})$ is the distance 
from $x$ to the boundary $\partial Q_{l,x}$ of $Q_{l,x}$. 
We remark for $x \in \Omega$,   $\tau_{0}(x)=\tau_{1}(x)$ if $\sup_{l \geq 0}
\displaystyle{\frac{\Delta_{l}(x)}{\Delta_{l+1}(x)}} < \infty$.

A following theorem may be proved by the same way as in [11].

\begin{Theorem}. 
Let $f$ and $\phi$ be bounded functions given in {\rm (18)}. 
Then we have

{\rm (i)} $\alpha (f,x) \geq \min (\alpha(\phi), \tau_{0}(x))$ for $x \in T$,

{\rm (ii)} $\alpha(f,x) \geq \min_{i}(\alpha(\phi, \Omega_{i}), \tau_{1}(x))$  
for $x \in \Omega$ with $\sup_{l \geq 0}
\displaystyle{\frac{\Delta_{l}(x)}{\Delta_{l+1}(x)}} < \infty$

\noindent
where $\Omega_{i} \equiv M^{-1}(T^{o}+ \gamma_{i})$, 
$\gamma_{i} \in \Gamma_{0}$, 
$i = 1, \ldots , m$ and $\alpha(\phi,\Omega_{i})=\sup \{ s\geq 0: \phi \in 
C^{s}(\Omega_{i}) \}$ and $C^{s}(\Omega_{i})$ is defined as the Besov 
space $B^{s}_{\infty \infty}(\Omega_{i})$ on $\Omega_{i}$.

{\rm (iii)} Suppose that $\phi \in C^{\infty}(\Omega_{i}), \ i=1, \ldots, m$ 
and  there exit a positive number $s_{0}$ and $y_{0} \in T^{o}$ such that 
\[
\sup_{l\geq 0}\sup_{y}\displaystyle{\frac{|f_{l}(y)|}
{(\lambda_{0}^{-l}+|y-y_{0}|)^{s_{0}}}}=\infty.
\]
 Then $\tau_{0}(x) \geq \alpha(f,x)$ for  $x \in T$.
\end{Theorem}

The proof of Theorem 4 is not difficult. We will omit its details.

\bigskip

\noindent
{\bf Corollary}. {\it
 Let $\phi$ be a bounded function on ${\Bbb R}^{n}$ such that 
$\phi \in C^{\infty}(\Omega_{j}), j=1,\ldots, m$ and $\phi=0$ outside $T^{o}$. 
Consider a bounded function $f$ given by {\rm (18)} satisfying the condition 
{\rm (iii)} in Theorem 4. Then we have

\noindent
{\rm (i)} $\ \ \ \tau_{0}(x) \geq \alpha(f,x) \geq 
\min(\alpha(\phi),\tau_{0}(x)) , \ \ x \in T$,

\noindent
{\rm (ii)}  for x in $\Omega$ with 
$\sup_{l \geq 0}\displaystyle{\frac{\Delta_{l}(x)}{\Delta_{l+1}(x)}}<\infty , 
\ \ \  
\alpha(f,x)= \tau_{0}(x)= \tau_{1}(x)$.  
}

\bigskip

\noindent 
{\bf Examples}.  
We consider a self-affine tiling $\Pi=\{ T+\nu \}_{\nu \in {\Bbb Z}}$ 
such that a tile $T=[0,1]$ and a dilation $M=2Id$ on ${\Bbb R}$.

\noindent
(a)  We consider the Takagi function such that 
\[
f(x)=\sum_{l=0}^{\infty}\sum_{Q \in \Pi_{l}(T)}\mu^{l}\phi(M_{Q}x) , 
\ \ \forall x \in {\Bbb R}
\]
where $0<\mu <1$ and $\phi$ is a bounded function such that 
$\phi(x)=x \ (0<x \leq \frac{1}{2}),\ \phi(x)=1-x \ (\frac{1}{2} \leq x <1),
\  \phi(x)=0 \ $(otherwise). 
Let $\tau=\displaystyle{\frac{\log \mu}{\log 2^{-1}}}$. Then 
from the corollary of Theorem 4, 
 if $\tau \leq 1$, $\tau=\alpha(f,x)$ 
for each $x \in T$.

\noindent
(b)  We consider the Weierstrass  function 
$f(x)=\sum_{l=0}^{\infty}\mu^{l}\phi(2^{l}x)$ with $0< \mu <1$ and 
$\phi(x)=\sin 2\pi x \ \ (x \in {\Bbb R})$.
 The proof of Theorem 4 can be also applied to this function 
  case. Then we have 
\[
\tau=\alpha(f,x), \ \ \forall x \in {\Bbb R}.
\]
where the constant $\tau=\displaystyle{\frac{\log \mu}{\log 2^{-1}}}$ 
is given in the part (a) above.

(c)  We consider  L${\grave {\rm e}}$vy's function 
\[
f(x)=\sum_{l=0}^{\infty}\sum_{Q \in \Pi_{l}(T)}2^{-l}\phi(M_{Q}x) , 
\ \ \forall x \in {\Bbb R}
\]
where $\phi(x)=x-\displaystyle{\frac{1}{2}}\ \ (0<x<1),\ \ 
\phi(x)=0 \ \  $(otherwise). 
 Then we can see  that       
$1=\tau_{1}(x)=\alpha(f,x)$ 
for a point $x$ in $\Omega$ 
 with $\sup_{l\geq 0} \displaystyle{\frac{\Delta_{l}(x)}{\Delta_{l+1}(x)}} 
 <\infty$.

{\sc Koichi Saka} \\
Department of Mathematics \\
Akita University \\
Akita, 010-8502 Japan \\
e-mail address: 
\verb+saka@math.akita-u.ac.jp+

\end{document}